\documentclass[12pt,reqno]{amsart}

\usepackage{mathrsfs}
\usepackage{amsmath}
\usepackage{amsfonts}
\usepackage{amssymb}
\usepackage[all,cmtip]{xy}
\usepackage{hyperref}
\usepackage{enumerate}  
\usepackage{enumitem}
\usepackage{ulem}
\usepackage[all]{xy}
\usepackage[toc,page]{appendix}
\usepackage{tikz-cd}
\usepackage[all,cmtip]{xy}
\usepackage[toc,page]{appendix}
\usepackage{amsmath, amsthm, amsfonts, amssymb}
\usepackage{graphicx}
\usepackage{subcaption}
\usepackage{tikz-cd}
\usepackage[utf8]{inputenc}
\usepackage[english]{babel}
\usepackage{dsfont}
\usepackage{xcolor}
\usepackage{arydshln}
\usepackage{faktor}
\usepackage{mathtools}
\usepackage{enumitem}
\usepackage{scalerel}
\usepackage{booktabs}

\def\endpf{\hbox{\vrule height1.5ex width.5em}}

\usepackage{dsfont}
\newcommand\xdashmapsto[2][]{\mathrel{\mapstochar\xdashrightarrow[#1]{#2}}}

\makeatletter
\newcommand*{\da@rightarrow}{\mathchar"0\hexnumber@\symAMSa 4B }
\newcommand*{\da@leftarrow}{\mathchar"0\hexnumber@\symAMSa 4C }
\newcommand*{\xdashrightarrow}[2][]{%
	\mathrel{%
		\mathpalette{\da@xarrow{#1}{#2}{}\da@rightarrow{\,}{}}{}%
	}%
}
\newcommand{\xdashleftarrow}[2][]{%
	\mathrel{%
		\mathpalette{\da@xarrow{#1}{#2}\da@leftarrow{}{}{\,}}{}%
	}%
}
\newcommand{\xdashdownarrow}[2][]{%
	\mathrel{%
		\mathpalette{\da@xarrow{#1}{#2}\da@downarrow{}{}{\,}}{}%
	}%
}
\newcommand*{\da@xarrow}[7]{%
	\sbox0{$\ifx#7\scriptstyle\scriptscriptstyle\else\scriptstyle\fi#5#1#6\m@th$}%
	\sbox2{$\ifx#7\scriptstyle\scriptscriptstyle\else\scriptstyle\fi#5#2#6\m@th$}%
	\sbox4{$#7\dabar@\m@th$}%
	\dimen@=\wd0 %
	\ifdim\wd2 >\dimen@
	\dimen@=\wd2 %
	\fi
	\count@=2 %
	\def\da@bars{\dabar@\dabar@}%
	\@whiledim\count@\wd4<\dimen@\do{%
		\advance\count@\@ne
		\expandafter\def\expandafter\da@bars\expandafter{%
			\da@bars
			\dabar@ 
		}%
	}%
	\mathrel{#3}%
	\mathrel{%
		\mathop{\da@bars}\limits
		\ifx\\#1\\%
		\else
		_{\copy0}%
		\fi
		\ifx\\#2\\%
		\else
		^{\copy2}%
		\fi
	}%
	\mathrel{#4}%
}
\makeatother

\usepackage{graphicx}

\newcommand\tikzmark[1]{%
  \tikz[remember picture,overlay]\coordinate (#1);}

\newcommand{\underbracedmatrixll}[2]{%
  \left(\;\hspace{-.27in}
  \smash[b]{\underbrace{
    \begin{matrix}#1\end{matrix}
  }_{#2}}
  \;\right.
  \vphantom{\underbrace{\begin{matrix}#1\end{matrix}}_{#2}}
}

\newcommand{\underbracedmatrixrr}[2]{%
  \left. \;
  \smash[b]{\underbrace{
    \begin{matrix}#1\end{matrix}
  }_{#2}}
  \;\hspace{-.32in}\right)
  \vphantom{\underbrace{\begin{matrix}#1\end{matrix}}_{#2}}
}

\newtheorem{theorem}{Theorem}[section]
\newtheorem{lemma}[theorem]{Lemma}
\newtheorem{corollary}[theorem]{Corollary}
\newtheorem{proposition}[theorem]{Proposition}

\theoremstyle{definition}
\newtheorem{question}[theorem]{Question}
\newtheorem{definition}[theorem]{Definition}
\newtheorem{example}[theorem]{Example}
\newtheorem{remark}[theorem]{Remark}

\newtheorem{definitionremark}[theorem]{Definition-Remark}
\date{}

\makeatletter
\let\@wraptoccontribs\wraptoccontribs
\makeatother

\begin{document}

\title[A vanishing theorem for the canonical blow-ups of  Grassmann manifolds]{\bf  A vanishing theorem for the canonical blow-ups of  Grassmann manifolds}

\author[Hanlong Fang and Songhao Zhu]{Hanlong Fang and Songhao Zhu}

\address{School of Mathematical Sciences, Peking University, Beijing 100871, China. }
\email{hlfang$ @$pku.edu.cn}
\address{Department of Mathematics, Rutgers University, Piscataway, NJ 08857, USA. }
\email{sz446$ @$math.rutgers.edu}

\vspace{3cm} \maketitle

\begin{abstract}
Let $\mathcal T_{s,p,n}$ be the canonical blow-up of the Grassmann manifold $G(p,n)$ constructed by blowing up the Pl\"ucker coordinate subspaces associated with the parameter $s$.  We prove  that the higher cohomology groups of the tangent bundle of $\mathcal T_{s,p,n}$ vanish. As an application, $\mathcal T_{s,p,n}$ is locally rigid in the sense of Kodaira-Spencer.
\end{abstract}

\tableofcontents

\section{Introduction}
Sheaf cohomology of vector bundles is a fundamental object studied in complex geometry and algebraic geometry. Based on such data, one may derive interesting geometric properties of the base varieties. For instance, geometers exploit appropriate vanishing theorems of higher cohomology  to construct global sections of vector bundles, and extend sections of vector bundles from subvarieties to the ambient spaces.  Kodaira-Spencer  theory relates the local deformation of the complex structure of a complex manifold $X$ to the first cohomology of its tangent bundle $H^1(X, T_X)$. Kuranishi proved that an element $\theta\in H^1(X, T_X)$ represents a local deformation if and only if its obstruction $[\theta,\theta]\in H^2(X, T_X)$ vanishes, and thus parametrized the local deformation of complex structures by the so called Kuranishi family. 

The study of cohomology groups of equivariant vector bundles on homogeneous manifolds has a long history dating back at least to the celebrated Borel-Weil-Bott theorem in 1950s, which gives explicit formulas in terms of the representations of the groups acting on the manifolds. Since then, various work has been done in extending the Borel-Weil-Bott theorem under different circumstances. An important direction of further generalization is to compute the cohomology  of vector bundles on a larger class of manifolds, that is, the spherical varieties.  Kato (\cite{Kat}) and Tchoudjem  (\cite{Tc}) settled the line bundle case for certain special spherical varieties (wonderful varieties in the sense of De consini-Procesi \cite{DP})   in terms of the weights of the corresponding lie algebras. 

Our paper stems from a systematic study of the canonical blow-ups of Grassmann manifolds (\cite{F}). It is an interesting family of  spherical varieties, which generalizes the notion of wonderful varieties to  homeward varieties. Recall that, by a result of Bott (\cite{Bo}),   every smooth homogeneous algebraic manifold $X$ over $\mathbb C$ has trivial higher cohomology groups of its tangent bundle. This implies that $X$ is locally rigid, or equivalently, any deformation $X_t$ parametrized by a complex
manifold $T$ with $X_0$ analytically
isomorphic to $X$, is holomorphically trivial.
Another application in the theory of $\mathcal D$-modules is that every regular
function on the cotangent bundle of $X$ is the symbol of a differential operator on
$X$ with regular coefficients (see \cite{BB}). Bien-Brion (Proposition 4.2 in \cite{BB}) generalized Bott's theorem to regular spherical  Fano manifolds.

We are interested in the following question, which naturally relates to the computation of the cohomology groups of the tangent bundle of $\mathcal T_{s,p,n}$.

\begin{question}
What is the Kuranishi family of $\mathcal T_{s,p,n}$?
\end{question}

By a result of Sano (\cite{Sa}), the local deformation of $\mathcal T_{s,p,n}$ is unobstructed, or equivalently, the Kuranishi family of $\mathcal T_{s,p,n}$ is smooth. Unfortunately, Bien-Brion's brilliant argument does not apply here directly, for in general $\mathcal T_{s,p,n}$ is only weak Fano instead of Fano. More precisely, the difficulty comes from the fact that the restriction of a big and numerical effective line bundle may fail to be big. An example illustrating this is to blow up a point in $\mathbb {CP}^2$, and then restrict the pull-back line bundle of $\mathcal O_{\mathbb {CP}^2}(1)$ to the exceptional divisor. 

We state our main result as follows.
\begin{theorem}\label{bignef}
Let $\mathcal T_{s,p,n}$ be the canonical blow up of Grassmann manifolds. Then,
\begin{equation}
 H^i(\mathcal T_{s,p,n}, T_{\mathcal T_{s,p,n}})=0   \,,\,\,\,i>0.
\end{equation} In particular, $\mathcal T_{s,p,n}$ is locally rigid.
\end{theorem}

Noticing that $\mathcal T_{p,p,2p}$ is isomorphic to Kausz's (\cite{Ka}) modular compactifications of general linear groups over $\mathbb C$, we have that 
\begin{corollary}
Let $KGL_p$ be Kausz's modular compactification of the general linear group $GL(p,\mathbb C)$. The higher cohomology of the tangent bundles of $KGL_p$ vanishes. In particular,  $KGL_p$ is locally rigid. 
\end{corollary}

We now  briefly describe the main idea of the proof. Notice that the argument in \cite{BB} used the ampleness of the anticanonical bundle only when applying the Kodaira vanishing theorem. Hence, it is natural to expect a finer result if one can replace the Kodaira vanishing theorem by the Kawamata-Viehweg vanishing theorem.  To deal with the difficulty that the restriction of a big line bundle fails to be big, we use the Van der Waerden representation (see \cite{F}) to extract the very explicit geometry of  $\mathcal T_{s,p,n}$. The crucial step is Lemma \ref{dgb}, which shows that the $B$-invariant divisors of the boundary divisors can be derived from the restriction of the $B$-invariant divisors of $\mathcal T_{s,p,n}$. Eventually,  computation  yields that in our case the restriction of the anticanonical bundle of  $\mathcal T_{s,p,n}$ to the components of the boundary divisors is indeed big and numerical effective, which is sufficient to apply the Kawamata-Viehweg vanishing theorem.

The organization of the paper is as follows. In \S 2,  we recall the construction of the canonical blow-ups of Grassmann manifolds and the basic properties following \cite{F}. In \S 3.1, we study the cone of  effective divisors of the components of the boundary divisor of $\mathcal T_{s,p,n}$. In \S 3.2, we first establish some numerical formulas of the restriction for the anticanonical bundles (the proof for the case $p=n-s$ or $s$ is left to Appendices \ref{dvanderl} and \ref{dvanderld}); then prove that the restriction of the anticanonical bundle of $\mathcal T_{s,p,n}$ to the components of its boundary divisor is big and numerical effective. Finally, we prove the main theorem in \S 3.3.  

For the reader's convenience, we recall in Appendix \ref{section:nc} the construction of the local coordinate charts used in this paper (the Van der Waerden representation) as well as an example illustrating this. In Appendices \ref{dvanderl} and \ref{dvanderld}, we provide a detailed proof of the numerical formulas for the restriction of the anticanonical bundles, when $p=n-s$ or $s$. 
\medskip

{\noindent\bf Acknowledgement.}  The first author would like to thank  Xin Fu, Zhan Li, and Jie Liu for helpful discussions.

\section{Basic construction and properties}
In this section, we will recall some notions and results in \cite{F}.
\subsection{Construction of \texorpdfstring{$\mathcal T_{s,p,n}$}{rr}}
Let $G(p,n)$, $0<p<n$, be the Grassmann manifold  consisting of complex $p$-planes in the complex $n$-space. 
Each point $x\in G(p,n)$ one to one corresponds to an equivalence class of $p\times n$ non-degenerate matrices, where the equivalence relation is induced by the matrix multiplication from the left by the elements of the general linear group $GL(p,\mathbb C)$. A matrix representative  $\widetilde x$ of $x$ is a matrix in the corresponding equivalence class.

Define an index set $\mathbb I_{p,n}$ by
\begin{equation}\label{ipq}
\mathbb I_{p,n}:=\big\{(i_1,\cdots,i_p)\in\mathbb Z^p\big|1\leq i_p<i_{p-1}<\cdots <\cdots<i_1\leq n\}.
\end{equation}
Denote by $[\cdots ,z_I,\cdots]_{I\in\mathbb I_{p,n}}$ the homogeneous coordinates for the complex projective space $\mathbb {CP}^{N_{p,n}}$ where $N_{p,n}=\frac{n!}{p! (n-p)!}-1$. For each index $I=(i_1,\cdots,i_p)\in\mathbb I_{p,n}$  and  a matrix representative  $\widetilde x$ of  $x\in G(p,n)$, denote by  $P_I(\widetilde x)$  the determinant of the submatrix of $\widetilde x$ consisting of the $i_1^{th}, \cdots,i_p^{th}$ 
columns. The Pl\"ucker embedding of  $G(p,n)$ into $\mathbb {CP}^{N_{p,n}}$  can be given  by
\begin{equation}\label{plucker}
\begin{split}
\mathfrak e:G(p,n)&\longrightarrow\mathbb {CP}^{N_{p,n}}\\
x&\mapsto [\cdots ,P_I(\widetilde x), \cdots]_{I\in\mathbb I_{p,n}}
\end{split}.
\end{equation}

For $0<s<n$ and $0\leq k\leq p$, define index sets  $\mathbb I_{s,p,n}^{k}$ by
\begin{equation}\label{ikk}
\mathbb I_{s,p,n}^{k}:=\big\{(i_1,\cdots,i_p)\in\mathbb Z^p\big|{\footnotesize 1\leq i_p<\cdots<i_{k+1}\leq s\,;s+1\leq i_{k}<i_{k-1}<\cdots<i_1\leq n}\}.
\end{equation}
We have a partition
\begin{equation}
\mathbb I_{p,n}=\bigsqcup_{k=0}^p\mathbb I_{s,p,n}^{k}\,\,.
\end{equation}
Consider linear subspaces of $\mathbb {CP}^{N_{p,n}}$ as follows.
\begin{equation}\label{subl}
\mathbb {CP}^{N^k_{s,p,n}}:=\big\{[\cdots ,z_I,\cdots]_{I\in\mathbb I_{p,n}}\in\mathbb {CP}^{N_{p,n}}\big|z_I=0\,,\,\,\forall I\notin\mathbb I_{s,p,n}^k \big\}\,,\,\,0\leq k\leq p\,,
\end{equation}
where $N^k_{s,p,n}$ is the cardinal number of the set $\mathbb I_{s,p,n}^{k}$ minus $1$; by a slight abuse of notation, we denote the corresponding homogeneous coordinates by  $[\cdots ,z_I,\cdots]_{I\in\mathbb I^k_{s,p,n}}$. Recall the following projection (rational) map  $F_s^k$ by dropping the coordinates whose indices are not in $\mathbb I_{s,p,n}^k$.
\begin{equation}\label{FK}
\begin{split}
F_s^k:\mathbb {CP}^{N_{p,n}}&\dashrightarrow\mathbb {CP}^{N^k_{s,p,n}} \\
[\cdots ,z_I,\cdots]_{I\in\mathbb I_{p,n}}&\xdashmapsto{}[\cdots, z_I,\cdots]_{I\in\mathbb I^k_{s,p,n}}.
\end{split}
\end{equation}
We make the convention that  $\mathbb{CP}^{N^k_{s,p,n}}$ is a point and $F_s^k$ is the trivial map when $N^k_{s,p,n}=0,-1$. We can thus define a rational map $\mathcal K_{s,p,n}:G(p,n)\dashrightarrow \mathbb {CP}^{N_{p,n}}\times\mathbb {CP}^{N^0_{s,p,n}}\times\cdots\times\mathbb {CP}^{N^p_{s,p,n}}$ by \begin{equation}\label{jl}
   \mathcal K_{s,p,n}:=(\mathfrak e, F_s^0\circ\mathfrak e,\cdots,F_s^p\circ\mathfrak e) \,;
\end{equation}
or equivalently, \begin{equation}\label{bpp}
K_{s,p,n}(x)=\big( [\cdots ,P_I(\widetilde x),\cdots]_{I\in\mathbb I_{p,n}},[\cdots,P_I(\widetilde x),\cdots]_{I\in\mathbb I^0_{s,p,n}},\cdots,[\cdots,P_I(\widetilde x),\cdots]_{I\in\mathbb I^p_{s,p,n}}\big)\,.
\end{equation}

\begin{definition}\label{ypt}
Assume that $0<p<n$ and $0<s<n$. Let $\mathcal T_{s,p,n}$ be the scheme-theoretic closure  of the birational image of $G(p,n)$ under  $\mathcal K_{s,p,n}$ in $\mathbb {CP}^{N_{p,n}}\times\mathbb {CP}^{N^0_{s,p,n}}\times\cdots\times\mathbb {CP}^{N^p_{s,p,n}}$. We call $\mathcal T_{s,p,n}$ the {\it canonical blow-up of $G(p,n)$} with respect to the parameter $s$. 
\end{definition}

\begin{example}[\cite{De}]\label{de}
Denote by $[x_1,\cdots, x_s,y_1,\cdots,y_{n-s}]$ the homogeneous coordinates for the projective space $\mathbb{CP}^{n-1}$. Then $\mathcal T_{s,1,n}$ is the blow-up of $\mathbb{CP}^{n-1}$ along the union of the disjoint linear subspaces $\mathbb{CP}^{s-1}$ and $\mathbb{CP}^{n-s-1}$; $\mathcal K_{s,1,n}$ is given by 
\begin{equation}
\mathcal K_{s,1,n}([x_1,\cdots, x_s,y_1,\cdots,y_{n-s}])=[x_1,\cdots, x_s,y_1,\cdots,y_{n-s}]\times [x_1,\cdots, x_s]\times[y_1,\cdots,y_{n-s}].
\end{equation}
$\mathcal T_{1,1,3}$ is the Hirzebruch surface $\Sigma_1$.
\end{example}
\begin{example}
$\mathcal T_{p,p,2p}$ is the modular compactification of the reductive group $GL(p,\mathbb C)$ constructed  by Kausz (\cite{Ka}).  
\end{example}

Take a subgroup $GL(s,\mathbb C)\times GL(n-s,\mathbb C)$  of $GL(n,\mathbb C)$ as follows.
\begin{equation}\label{sns}
GL(s,\mathbb C)\times GL(n-s,\mathbb C):=\left\{\left.\left(
\begin{matrix}
g_1&0\\
0&g_2\\
\end{matrix}\right)\right\vert_{} g_1\in GL(s,\mathbb C), g_2\in GL(n-s,\mathbb C)\right\}.
\end{equation}
Let $B$ be a Borel subgroup of $GL(s,\mathbb C)\times GL(n-s,\mathbb C)$ given by \begin{equation}
B:=\left\{\left. \left(
\begin{matrix}
g_1&0\\
0&g_2\\
\end{matrix}\right)\right\vert_{} 
\begin{matrix}
g_1\in GL(s,\mathbb C) \,\, {\rm is\,\, a\,\, lower\,\, triangular\,\, matrix\,};\,\,\,\,\\
g_2\in GL(n-s,\mathbb C)\,\, {\rm is\,\, an \,\, upper \,\, triangular \,\, matrix}  
\end{matrix}
\right\}.
\end{equation} 

\begin{definition} Let $G$ be a connected reductive group. An irreducible normal $G$-variety $X$ is called {\it spherical} if a Borel subgroup of $G$ has an open orbit on $X$. \end{definition}

Recall that
\begin{proposition}[Propositions 1.3 and 1.12 in \cite{F}]\label{gwond}
$\mathcal T_{s,p,n}$ is a smooth spherical $GL(s,\mathbb C)\times GL(n-s,\mathbb C)$-variety. The complement of the open
$GL(s,\mathbb C)\times GL(n-s,\mathbb C)$-orbit in $\mathcal T_{s,p,n}$ is a simple normal crossing divisor consisting of $2r$  
smooth, irreducible divisors as follows.
\begin{equation}    
D^-_1, 
D^-_2,\cdots,D^-_r,D^+_1,D^+_2,\cdots,D^+_r\,.
\end{equation}
Each $GL(s,\mathbb C)\times GL(n-s,\mathbb C)$-orbit of $\mathcal T_{s,p,n}$ one to one corresponds to the  quasi-projective variety $X_{(I^-,I^+)}$ defined  by
\begin{equation}\label{inrule}
X_{(I^-,I^+)}:=\left(\bigcap_{\,\,i\in I^-}D^-_i\mathbin{\scaleobj{1.2}{\bigcap}} \bigcap_{\,\,i\in I^+}D^+_i\right)\mathbin{\scaleobj{2.1}{\backslash}} \left(\bigcup_{\substack{1\leq j\leq r\\\,\,j\notin I^-}}D^-_j\mathbin{\scaleobj{1.2}{\bigcup}}\,\bigcup_{\substack{1\leq j\leq r\\\,\,\,j\notin I^+}}D^+_j\right),
\end{equation} 
where $I^-,I^+$ are subsets of $\{1,2,\cdots,r\}$
such that 
\begin{equation}
\min(I^-)+\min(I^+)\geq r+2.
\end{equation}
Here we make the convention that $\min(\emptyset)=+\infty$.
Moreover, the closure of each $GL(s,\mathbb C)\times GL(n-s,\mathbb C)$-orbit in $\mathcal T_{s,p,n}$ is smooth.
\end{proposition}

\begin{definitionremark}
One can show that \begin{equation}
    r=\min\{s,n-s,p,n-p\}. 
\end{equation} We make the convention that  $r$ is always referred to the above quantity in this paper.
\end{definitionremark}
A line bundle $L$
on a projective variety $X$ of dimension $n$ is called big if its highest self-intersection number $(L^n)$ is positive, and called numerical effective (or nef for short) if the intersection number $(L\cdot C)$ is non-negative for any complete curve $C$ on $X$.

We have that

\begin{proposition}[Theorem 1.22 in \cite{F}]\label{fano}
The anti-canonical bundle $-K_{\mathcal T_{s,p,n}}$ of $\mathcal T_{s,p,n}$ is big and numerical effective. $-K_{\mathcal T_{s,p,n}}$ is ample if and only if $r\leq 2$. 
\end{proposition}

According to Definition-Lemmas 2.11 and 2.12 in \cite{F}, the following isomorphisms among the canonical blow-ups of Grassmann manifolds hold. 
\begin{equation}
{\rm DUAL}:\,\mathcal T_{s,p,n}\rightarrow\mathcal T_{s,n-p,n}\,\,\,\,{\rm and}\,\,\,\,
     {\rm USD}:\mathcal T_{s,p,n}\rightarrow\mathcal T_{n-s,p,n}\,.\\
\end{equation}
\begin{remark}
Without loss of generality, in the remaining of this paper, we can thus assume that $2p\leq n\leq 2s$.
\end{remark}

\subsection{$G$-invariant divisors of \texorpdfstring{$\mathcal T_{s,p,n}$}{rr}}\label{basict}

In this subsection, we will give a more detailed description of the divisors $D^-_1, 
D^-_2,\cdots,D^-_r,D^+_1,D^+_2,\cdots,D^+_r$   appearing in Proposition \ref{gwond}. For convenience, we denote by
$G$ the group $GL(s,\mathbb C)\times GL(n-s,\mathbb C)$ in the following.

Define an algebraic $\mathbb C^*$-action $\psi_{s,p,n}$ on $G(p,n)$  by
\begin{equation}\label{grac} \psi_{s,p,n}(\lambda):=\left(\begin{matrix} I_{s\times s}&0\\ 0&\lambda\cdot I_{(n-s)\times (n-s)}\\ \end{matrix} \right),\,\,\lambda\in\mathbb C^*.
\end{equation}
We have a unique lifting $\Psi_{s,p,n}$ of $\psi_{s,p,n}$ from $G(p,n)$ to $\mathcal T_{s,p,n}$ (see Lemma 2.10 in \cite{F}).

For $0\leq l\leq r$, define subsets $\mathcal V_{(p-l,l)}$, $\mathcal V_{(p-l,l)}^+$ and $\mathcal V_{(p-l,l)}^-$ of $G(p,n)$ in matrix representatives by
\begin{equation}\label{vani}
\begin{split}
&\mathcal V_{(p-l,l)} :=\left\{\left. \left(
\begin{matrix}
0&X\\
Y&0\\
\end{matrix}\right)\right\vert_{}\footnotesize\begin{matrix}
X\,\,{\rm is\,\,an\,\,}l\times (n-s)\,\,{\rm matrix\,\,of\,\,rank}\,\,l\,;\\
Y\,\,{\rm is\,\,a\,\,}(p-l)\times s\,\,{\rm matrix\,\,of\,\,rank}\,\,(p-l)\\
\end{matrix}
\right\}\,,
\end{split}
\end{equation}
and
\begin{equation}\label{stab}
\begin{split}
&\mathcal V_{(p-l,l)}^+:= \left\{\left.\left(
\begin{matrix}
0&X\\
Y&W\\
\end{matrix}\right)\right\vert_{}\footnotesize\begin{matrix}
X\,\,{\rm is\,\,an\,\,}l\times (n-s)\,\,{\rm matrix \,\,of\,\,rank\,\,}l\,;\\
Y\,\,{\rm is\,\,a\,\,}(p-l)\times s\,\,{\rm matrix \,\,of\,\,rank\,\,}(p-l)\,\\
\end{matrix}\right\},\\
&\mathcal V_{(p-l,l)}^-:= \left\{\left.\left(
\begin{matrix}
Z&X\\
Y&0\\
\end{matrix}\right)\right\vert_{}{\footnotesize\begin{matrix}
X\,\,{\rm is\,\,an\,\,}l\times (n-s)\,\,{\rm matrix \,\,of\,\,rank\,\,}l\,;\\
Y\,\,{\rm is\,\,a\,\,}(p-l)\times s\,\,{\rm matrix \,\,of\,\,rank\,\,}(p-l)\,\\
\end{matrix}}\right\}\,.    
\end{split}
\end{equation}

We have the following explicit  Bia{\l}ynicki-Birula decomposition (see \cite{Bi}) for $\mathcal T_{s,p,n}$.

\begin{lemma}[Lemma 4.9 in \cite{F}]\label{11} 
There are $r+1$ connected components $\mathcal D_{(p-l,l)}$, $0\leq l\leq r$, of the  set of the fixed points of $\mathcal T_{s,p,n}$ under  the $\mathbb C^*$-action $\Psi_{s,p,n}$, such that  the following holds.
\begin{enumerate}[label=(\alph*).]
\item   $R_{s,p,n}(\mathcal D_{(p-l,l)})=\mathcal V_{(p-l,l)}$, $0\leq l\leq r$. 

\item For $0\leq l\leq r$, there is a fibration $\mathcal D^+_{(p-l,l)}$ (resp. $\mathcal D^-_{(p-l,l)}$) over $D_{(p-l,l)}$ such that  
\begin{equation}\label{tschubert}
\overline{\mathcal D_{(p-l,l)}^+}\mathbin{\scaleobj{1.7}{\backslash}} \mathcal D_{(p-l,l)}^+=\bigsqcup_{k=l+1}^r\mathcal D_{(p-k,k)}^+\,\,\,\,\,\,\,\,\left({\rm resp.}\,\,\,\,\,\,\,\overline{\mathcal D_{(p-l,l)}^-}\mathbin{\scaleobj{1.7}{\backslash}} \mathcal D_{(p-l,l)}^-=\bigsqcup_{k=0}^{l-1}\mathcal D_{(p-k,k)}^-\,\right)\,.
\end{equation}
	
\end{enumerate}
\end{lemma}

\begin{lemma}[see Definition 4.10 and the proof of Proposition 1.12 in \cite{F}]
For $1\leq k\leq r$,  $D_{k}^-$ (resp. $D_{k}^+$) is the Zariski closure of the manifold $\mathcal D_{(p-k+1,k-1)}^-$  $\left({\rm resp}.\,\,\mathcal D_{(p-r+k-1,r-k+1)}^+\right)$ .
\end{lemma}

\subsection{$B$-invariant divisors of \texorpdfstring{$\mathcal T_{s,p,n}$}{rr}}
To describe the cone of effective divisors of $\mathcal T_{s,p,n}$, we need to find its $B$-invariant divisors. In this subsection, we will recall some important properties of the $B$-invariant divisors of $\mathcal T_{s,p,n}$.

For $0\leq j\leq r$, define irreducible divisors $b_j$  of $G(p,n)$ by 
\begin{equation}\label{ij}
 b_{j}:=\left\{x\in G(p,n) \big|\,P_{I_j}(x)=0\,,\,\,I_j=(s+j,s+j-1,\cdots,s-p+j+1)\in\mathbb I^j_{s,p,n}\right\}\,.
\end{equation}
Let $B_{j}\subset\mathcal T_{s,p,n}$ be the strict transformation of $b_j$ under the blow-up $R_{s,p,n}:\mathcal T_{s,p,n}\rightarrow G(p,n)$.  Notice that when $p=n-s$, $B_r=D^-_r$; when $p=s$, $B_0=D^+_r$.

Recall that 
\begin{lemma}[Lemma 6.7 in \cite{F}]\label{bb}When $0\leq j\leq r$,
\begin{equation}\label{bst}
\begin{split}
 B_j&=(R_{s,p,n})^*\left(\mathcal O_{G(p,n)}(1)\right)-\sum_{i=1}^{r-j}(r-j+1-i)\cdot D^+_i-\sum_{i=1}^j(j+1-i)\cdot D^-_i.
\end{split}
\end{equation} If $p=s$ (resp. $p=n-s$) we should modify (\ref{bst}) for $B_0$ (resp. $B_r$) as follows. When $p=s$, 
\begin{equation}\label{bstt0}
\begin{split}
B_0=D^+_r&=(R_{s,p,n})^*\left(\mathcal O_{G(p,n)}(1)\right)-\sum_{i=1}^{r-1}(r+1-i)\cdot D^+_i\,;
\end{split}
\end{equation}
when $p=n-s$,
\begin{equation}\label{bsttr}
\begin{split}
  B_r=D^-_r&=(R_{s,p,n})^*\left(\mathcal O_{G(p,n)}(1)\right)-\sum_{i=1}^{r-1}(r+1-i)\cdot D^-_i\,.
\end{split}
\end{equation}
\end{lemma}

And
\begin{lemma}[Lemma 3.3 in \cite{F}]\label{orbits}
If $B_j$ contains a non-empty $G$-orbit of  $\mathcal T_{s,p,n}$, then  either  $p=n-s$, $j=r$ and $B_j=B_r=D^-_r$, or $p=s$, $j=0$ and  $B_j=B_0=D^+_r$. 
\end{lemma} 

Immediately, we have that 
\begin{lemma}\label{regu}
$\mathcal T_{s,p,n}$ is regular in the sense of \cite{BB}. Precisely, $\mathcal T_{s,p,n}$ is smooth and spherical  without color (i.e. every irreducible $B$-stable divisor containing a $G$-orbit is $G$-stable).
\end{lemma}

For a smooth projective manifold $X$ over $\mathbb C$,  the group $Z_{i}(X)$ of $i$-dimensional cycles on $X$ is the free abelian group on the set of $i$-dimensional subvarieties of $X$; the group of $i$-cycles rationally equivalent to zero is the subgroup of $Z_{i}(X)$ generated by the cycles $(f)$ for all $(i+1)$-dimensional subvarieties $W$ of $X$ and all nonzero rational functions $f$ on $W$; the Chow group $A_{i}(X)$ of $i$-dimensional cycles on $X$ is the quotient group of $Z_{i}(X)$ by the subgroup of cycles rationally equivalent to zero. 

Brion (\cite{Br}) proved that
\begin{lemma}\label{cone}
Let $X$ be an irreducible, complete spherical variety of complex dimension $n$.
The cone of effective divisors in $A_{n-1}(X)\otimes_{\mathbb Z}\mathbb Q$ is a polyhedral convex cone generated by the classes of irreducible $B$-invariant divisors.
\end{lemma}

We determine the  cone of effective divisors of $\mathcal T_{s,p,n}$ by 

\begin{lemma}[Lemma 3.2 in \cite{F}]\label{gb}
Let $\mathfrak D$ be an irreducible $B$-invariant divisor of $\mathcal T_{s,p,n}$. Then
\begin{equation}\label{bi}
\mathfrak D\in \{D_{1}^-, D_{2}^-, \cdots, D_{r}^-, D_{1}^+, D_{2}^+, \cdots, D_{r}^+, B_0, B_1, \cdots, B_r\}\,.
\end{equation}
\end{lemma}
\begin{remark}
Since $b_m$ is biholomporphic to the infinity hyperplane section of $G(p,n)$ which is the closure of a complex Euclidean space, $B_m$ is  irreducible and $B$-invariant.
\end{remark}
{\noindent\bf Proof of Lemma \ref{gb}.} For the readers' convenience, we repeat here the proof given in \cite{F}. 

Assume that $\mathfrak D$ is $B$-invariant but $\mathfrak D\not\in\left\{D_{1}^-,\cdots, D_{r}^-, D_{1}^+, \cdots, D_{r}^+, B_0,\cdots,B_r\right\}$.  Denote by $\mathfrak d$ the image of $\mathfrak D$ under $R_{s,p,n}$. It is clear that $\mathfrak d$ is a $B$-invariant divisor of $G(p,n)$, for the exceptional divisor of $R_{s,p,n}$ is contained in the union of the $G$-stable divisors. Hence, for $0\leq j\leq r$, $P_{I_j}\not\equiv 0$ on $\mathfrak d$  (see (\ref{ij}) for the definition). We can verify that $\mathfrak d$ contains a point $\mathfrak a$ with 
a matrix representative $\widetilde {\mathfrak a}$ defined by \begin{equation}\label{1pp}
\widetilde {\mathfrak a}:=\left(0_{p\times(s-p)}
\hspace{-.13in}\begin{matrix}
  &\hfill\tikzmark{g}\\
  &\hfill\tikzmark{h}
  \end{matrix}\,\,\,\,\begin{matrix}
  I_{p\times p} \hspace{-.1in}\begin{matrix}
  &\hfill\tikzmark{c}\\
  &\hfill\tikzmark{d}
  \end{matrix}\hspace{-.1in}\begin{matrix}
  &\hfill\tikzmark{e}\\
  &\hfill\tikzmark{f}\end{matrix}\,\,\,\,\,I_{p\times p}
  \end{matrix}\hspace{-.1in}
  \begin{matrix}
  &\hfill\tikzmark{a}\\
  &\hfill\tikzmark{b}\end{matrix}\,\,\,\,\,0_{p\times(n-s-p)}\right)\,\,{\rm when}\,\, r=p\leq n-s,
  \tikz[remember picture,overlay]   \draw[dashed,dash pattern={on 4pt off 2pt}] ([xshift=0.5\tabcolsep,yshift=7pt]a.north) -- ([xshift=0.5\tabcolsep,yshift=-2pt]b.south);\tikz[remember picture,overlay]   \draw[dashed,dash pattern={on 4pt off 2pt}] ([xshift=0.5\tabcolsep,yshift=7pt]c.north) -- ([xshift=0.5\tabcolsep,yshift=-2pt]d.south);\tikz[remember picture,overlay]   \draw[dashed,dash pattern={on 4pt off 2pt}] ([xshift=0.5\tabcolsep,yshift=7pt]e.north) -- ([xshift=0.5\tabcolsep,yshift=-2pt]f.south);\tikz[remember picture,overlay]   \draw[dashed,dash pattern={on 4pt off 2pt}] ([xshift=0.5\tabcolsep,yshift=7pt]g.north) -- ([xshift=0.5\tabcolsep,yshift=-2pt]h.south);
\end{equation}
or
\begin{equation}\label{1pr}
\widetilde {\mathfrak a}:=\left(\begin{matrix}
  0_{r\times(s-p)}\\
   0_{(p-r)\times(s-p)}
 \end{matrix}
\hspace{-.13in}\begin{matrix}
  &\hfill\tikzmark{g}\\
  &\hfill\tikzmark{h}
  \end{matrix}\,\,\,\,\begin{matrix}
  I_{r\times r}\\
   0_{(p-r)\times r}
  \end{matrix}\hspace{-.1in}\begin{matrix}
  &\hfill\tikzmark{c}\\
  &\hfill\tikzmark{d}
  \end{matrix}\,\,\,\,\,\begin{matrix}
  0_{r\times(p-r)}\\
   I_{(p-r)\times(p-r)}
  \end{matrix}\hspace{-.1in}
\begin{matrix}
  &\hfill\tikzmark{e}\\
  &\hfill\tikzmark{f}\end{matrix}\hspace{-.1in}
  \begin{matrix}
  &\hfill\tikzmark{a}\\
  &\hfill\tikzmark{b}\end{matrix}\,\,\,\,\,
  \begin{matrix}
  I_{r\times r}\\
   0_{(p-r)\times r}
  \end{matrix}\right)\,\,{\rm when}\,\,r=n-s\leq p. 
  \tikz[remember picture,overlay]   \draw[dashed,dash pattern={on 4pt off 2pt}] ([xshift=0.5\tabcolsep,yshift=7pt]a.north) -- ([xshift=0.5\tabcolsep,yshift=-2pt]b.south);\tikz[remember picture,overlay]   \draw[dashed,dash pattern={on 4pt off 2pt}] ([xshift=0.5\tabcolsep,yshift=7pt]c.north) -- ([xshift=0.5\tabcolsep,yshift=-2pt]d.south);\tikz[remember picture,overlay]   \draw[dashed,dash pattern={on 4pt off 2pt}] ([xshift=0.5\tabcolsep,yshift=7pt]e.north) -- ([xshift=0.5\tabcolsep,yshift=-2pt]f.south);\tikz[remember picture,overlay]   \draw[dashed,dash pattern={on 4pt off 2pt}] ([xshift=0.5\tabcolsep,yshift=7pt]g.north) -- ([xshift=0.5\tabcolsep,yshift=-2pt]h.south);
\end{equation}
Then $\mathfrak a$ is in a dense open $B$-orbit of $G(p,n)$, which is a contradiction.

We complete the proof of Lemma \ref{gb}.\,\,\,$\endpf$.

\section{Semi-positivity of the restriction of the anticanonical bundles}
Throughout this section, we assume that $r\geq 3$ and $2p\leq n\leq 2s$.

\subsection{Geometric structure of the  cone of effective divisors of \texorpdfstring{$D_j^{\pm}$}{jj}} 

For $1\leq j\leq r$ and $0\leq m\leq r$, denote by $\check B^{-j}_{m}$ (resp. $\check B^{+j}_{m}$) the restriction of the line bundle $B_m$ to $D^{-}_j$ (resp. $D^{+}_j$). Notice that when $B_m\neq D^{-}_j$ (resp. $B_m\neq D^{+}_j$), we can identify $\check B^{-j}_{m}$ (resp. $\check B^{+j}_{m}$) with an effective divisor of $D^{-}_j$ (resp. $D^{+}_j$), that is the scheme-theoretic intersection of $B_m$ and $D^{-}_j$ (resp. $D^{+}_j$).  For $1\leq i, j\leq r$, denote by $D^{\pm j}_{\pm i}$ the restriction of the line bundle $D^{\pm}_i$ to $D^{\pm}_j$. 
Similarly, when $D^{\pm}_i\neq D^{\pm}_j$, we can identify $D^{\pm j}_{\pm i}$ with an effective divisor of $D^{\pm}_j$, that is, the scheme-theoretic intersection of $D^{\pm}_j$ and $D^{\pm}_i$. We make the convention that when the indices are out of the above range,  $D^{\pm i}_{\pm j}$ and $\check B^{\pm j}_{m}$ represent the trivial line bundles.

We prove the following crucial lemma similarly to Lemma \ref{gb}. 

\begin{lemma}\label{dgb}

Let $\mathfrak D$ be an irreducible divisor of $D^{-}_j$ (resp. $D^{+}_j$), $1\leq j\leq r$. If $\mathfrak D$ is $G$-invariant, then
\begin{equation}\label{dgi}
\begin{split}
&\,\,\,\,\,\,\,\,\,\,\mathfrak D\in \left\{D^{- j}_{- 1}, D^{- j}_{- 2}, \cdots,\widehat{D^{- j}_{- j}}, \cdots, D^{- j}_{- r}, D^{- j}_{+ (r+2-j)}, D^{- j}_{+ (r+3-j)}, \cdots, D^{- j}_{+ r}\right\}\,\\
&\left( resp.\,\, \mathfrak D\in \left\{D^{+ j}_{+ 1}, D^{+ j}_{+ 2}, \cdots,\widehat{D^{+ j}_{+ j}}, \cdots, D^{+ j}_{+ r}, D^{+ j}_{- (r+2-j)}, D^{+ j}_{- (r+3-j)}, \cdots, D^{+ j}_{- r}\right\}\right)\,.\\
\end{split}
\end{equation}
If $\mathfrak D$ is $B$-invariant but not $G$-invariant, then for $p\neq n-s$ and $p\neq s$, $\mathfrak D$ is an irreducible component of one of the following divisors,
\begin{equation}\label{dbi1}
 \left\{\check B^{- j}_{0}, \check B^{- j}_{1}, \cdots, \check B^{- j}_{r}\right\}\,\,\,\,\left(resp.\,\,\left\{\check B^{+ j}_{0}, \check B^{+ j}_{1}, \cdots, \check B^{+ j}_{r}\right\}\right);
\end{equation}
for $p=n-s<s$, $\mathfrak D$ is an irreducible component of one of the following divisors,
\begin{equation}\label{dbi2}
\left\{\check B^{- j}_{0}, \check B^{- j}_{1}, \cdots, \check B^{- j}_{r-1}\right\}\,\,\,\,\left(resp.\,\,\left\{\check B^{+ j}_{0}, \check B^{+ j}_{1}, \cdots, \check B^{+ j}_{r-1}\right\}\right);
\end{equation}
for $p=n-s=s$,  $\mathfrak D$ is an irreducible component of one of the following divisors,
\begin{equation}\label{dbi3}
\left\{\check B^{- j}_{1}, \check B^{- j}_{2}, \cdots, \check B^{- j}_{r-1}\right\}\,\,\,\,\left(resp.\,\,\left\{\check B^{+ j}_{1}, \check B^{+ j}_{2}, \cdots, \check B^{+ j}_{r-1}\right\}\right).
\end{equation}
\end{lemma}


{\bf\noindent Proof of Lemma
\ref{dgb}.} (\ref{dgi}) follows from Proposition \ref{gwond}.

The remaining of the proof is similar to that of Lemma \ref{gb}.

Firstly, assume that $p\neq n-s$ or $s$. Without loss of generality, we can consider divisors of  $D^{-}_j$, $1\leq j\leq r$; the case of $D^{+}_j$ is similar, and we omit it here for simplicity.

Recall that $D^{-}_1$,  $D^{-1}_m$ and $\check B^{-1}_m$ are $\mathcal M_{s,p,n}$, $\check D_{m}$, and $\check B_m$ defined in \cite{F}, respectively. (\ref{dbi1}) is exactly  Lemma 6.16 in \cite{F}. 

In what follows, we will assume  $2\leq j\leq r$ and proceed to prove (\ref{dbi1}) by contradiction. Suppose that ${\mathfrak D}$ is $B$-invariant but not $G$-invariant, and $\mathfrak D\not\in\{\check B^{-j}_{0}, \check B^{- j}_{1}, \cdots, \check B^{-j}_{r}\}$. Let ${\mathfrak a}$ be a generic point of ${\mathfrak D}$; that is, ${\mathfrak a}$ is not a point of $D^{- j}_{- 1}, D^{- j}_{- 2}, \cdots,\widehat{D^{- j}_{- j}}, \cdots, D^{- j}_{- r}, D^{- j}_{+ (r+2-j)}, D^{- j}_{+ (r+3-j)}$, $\cdots, D^{- j}_{+ r},\check B^{-j}_{0}, \check B^{- j}_{1}, \cdots, \check B^{-j}_{r}$. Denote by $\bar{\mathfrak a}$ the image of $\mathfrak a$ under $R_{s,p,n}:\mathcal T_{s,p,n}\rightarrow G(p,n)$. Then $\bar{\mathfrak a}$ has the following matrix representative,
\begin{equation}
\widetilde{\bar{\mathfrak a}}=\left(
\begin{matrix}
Z&X\\
Y&0\\
\end{matrix}\right)\,,\,\,\,
\end{equation}
where $X$ is an $(j-1)\times(n-s)$ matrix of rank $(j-1)$ and $Y$ is a $(p-j+1)\times s
$ matrix of rank $(p-j+1)$. Recall (\ref{ij}) that \begin{equation}
 b_{j-1}=\left\{x\in G(p,n) \big|\,P_{I_{j-1}}(x)=0\,,\,\,I_j=(s+j-1,s+j-2,\cdots,s-p+j)\right\}\,;
\end{equation}
by (\ref{bst}), the preimage of $b_{j-1}$ under the map $R_{s,p,n}$ is the union of $B_{j-1}$, $D^-_{-1},D^-_{-2},\cdots$, $D^-_{j-1}$, $D^+_1$, $D^+_2,\cdots, D^+_{r-j+1}$. Therefore, $\widetilde{\bar{\mathfrak a}}\notin b_{j-1}$,  $P_{I_{j-1}}(\widetilde{\bar{\mathfrak a}})\neq 0$, and we can conclude that $\bar{\mathfrak a}$ has a matrix representative as follows. 
\begin{equation}
\widetilde{\bar{\mathfrak a}}=\underbracedmatrixll{Z\\Y}{s-p+j-1\,\,\rm columns}
  \hspace{-.45in}\begin{matrix}
  &\hfill\tikzmark{a}\\
  &\hfill\tikzmark{b}  
  \end{matrix} \,\,\,\,\,
  \begin{matrix}
  0\\
I_{(p-j+1)\times(p-j+1)}\\
\end{matrix}\hspace{-.11in}
\begin{matrix}
  &\hfill\tikzmark{c}\\
  &\hfill\tikzmark{d}
  \end{matrix}\hspace{-.11in}\begin{matrix}
  &\hfill\tikzmark{g}\\
  &\hfill\tikzmark{h}
  \end{matrix}\,\,\,\,
\begin{matrix}
I_{(j-1)\times (j-1)}\\
0\\
\end{matrix}\hspace{-.11in}
\begin{matrix}
  &\hfill\tikzmark{e}\\
  &\hfill\tikzmark{f}\end{matrix}\hspace{-.3in}\underbracedmatrixrr{X\\0}{(n-s-j+1)\,\,\rm columns}\,\,\,,
  \tikz[remember picture,overlay]   \draw[dashed,dash pattern={on 4pt off 2pt}] ([xshift=0.5\tabcolsep,yshift=7pt]a.north) -- ([xshift=0.5\tabcolsep,yshift=-2pt]b.south);\tikz[remember picture,overlay]   \draw[dashed,dash pattern={on 4pt off 2pt}] ([xshift=0.5\tabcolsep,yshift=7pt]c.north) -- ([xshift=0.5\tabcolsep,yshift=-2pt]d.south);\tikz[remember picture,overlay]   \draw[dashed,dash pattern={on 4pt off 2pt}] ([xshift=0.5\tabcolsep,yshift=7pt]e.north) -- ([xshift=0.5\tabcolsep,yshift=-2pt]f.south);\tikz[remember picture,overlay]   \draw[dashed,dash pattern={on 4pt off 2pt}] ([xshift=0.5\tabcolsep,yshift=7pt]g.north) -- ([xshift=0.5\tabcolsep,yshift=-2pt]h.south);
\end{equation}
where
\begin{equation}
\small
    Z:=\left(\begin{matrix}
    z_{11}&\cdots&z_{1(s-p+j-1)}\\
    \vdots&\ddots&\vdots\\
    z_{l1}&\cdots&z_{l(s-p+j-1)}\\
    \end{matrix}\right),X:=\left(\begin{matrix}x_{1(s+j)}&\cdots&x_{1n}\\ \vdots&\ddots&\vdots\\ x_{l(s+j)}&\cdots &x_{ln}\\
\end{matrix}\right),Y:=\left(\begin{matrix}
y_{j1}&\cdots& y_{j(s-p+j-1)}\\
\vdots&\ddots&\vdots\\
y_{p1}&\cdots& y_{p(s-p+j-1)}\\
\end{matrix}  \right).
\end{equation}

By the Van der Waerden representation presented in \cite{F} (see Appendiex \ref{section:nc} as well), we can  take a holomorphic coordinate chart $\left(A^{\tau},(J_{j-1}^{\tau})^{-1}\right)$ for $\tau\in\mathbb J_{j-1}$ (see (\ref{tt}) for the definition of $\mathbb J_{j-1}$) around $\mathfrak a$ in $\mathcal T_{s,p,n}$,
such that the projection $R_{s,p,n}$ in the local coordinates  $\left(\widetilde X,\widetilde Y,\overrightarrow B^1,\cdots,\overrightarrow B^p\right)$ (see (\ref{ulu}), (\ref{qbp}) and (\ref{qbpb})) takes the following form.
\begin{equation}
\begin{split}
&\,\,\,\,\,\,\,\,\,\,\,\,\,\,\,\,\,\,\,\,\,\,\,\,\,\,\,\,\,\,\,\,\,\,\,\,\,\,\,\,\,\,\,\,\,\,\,\,\,\,\,\,\,\,\,\,\,\,\,\,\,\,\,\,\,\Gamma_{j-1}^{\tau}\left(\widetilde X,\widetilde Y,\overrightarrow B^1,\cdots,\overrightarrow B^p\right):=\\
&\small\left(
\begin{matrix}
\sum\limits_{k=p-j+2}^p\left(\prod\limits_{t=p-j+2}^{k}a_{i_{t}j_t}\right)\cdot\Xi_k^T\cdot\Omega_k &0_{(j-1)\times(p-j+1)}&I_{(j-1)\times (j-1)}&\widetilde X\\ \widetilde Y&I_{(p-j+1)\times(p-j+1)}&0_{(p-j+1)\times (j-1)}&\sum\limits_{k=1}^{p-j+2}\left(\prod\limits_{t=1}^{k}b_{i_tj_t}\right)\cdot\Xi_k^T\cdot\Omega_k\\
\end{matrix}\right)\,,  
\end{split}
\end{equation}
where $\Xi_k$ and $\Omega_k$ are given by (\ref{w3}), (\ref{w4}), (\ref{w5}) and (\ref{w2}). 
\smallskip

{\noindent\bf Claim.} We can choose the index $\tau\in\mathbb J_{j-1}$  in the above as 
\begin{equation}
\tau=\left(\begin{matrix}
j&j+1&\cdots&p&j-1&j-2&\cdots&1\\
s+j&s+j+1&\cdots&s+p&s-p+j-1&s-p+j-2&\cdots& s-p+1\\
\end{matrix}\right).
\end{equation}

{\noindent\bf Proof of Claim.}
Following \cite{F}, we define special indices $I_k,I_k^*,  I_{\mu\nu}^k, I_{\mu\nu}^{k*}\in\mathbb I^k_{s,p,n}$ as follows. For $0\leq k\leq p$, define
\begin{equation}\label{I_k}
    I_k:=(s+k,s+k-1,\cdots,s-p+k+1)\,;
\end{equation}
for $1\leq k\leq p-1$, define
\begin{equation}\label{I*k}
    I^*_{k}:=(s+k+1,s+k-1,s+k-2,\cdots,s-p+k+3,s-p+k+2,s-p+k)\,;
\end{equation}
for $0\leq k\leq p$, $s-p+k+1\leq \mu\leq s$, and $1\leq \nu\leq s-p+k$ define
\begin{equation}\label{Ikmn}
    I_{\mu\nu}^k:=(s+k,s+k-1,\cdots,,\widehat \mu,\cdots, \nu)\,;
\end{equation}
for $0\leq k\leq p$, $s+1\leq \mu\leq s+k$, and $s+k+1\leq \nu\leq n$ define
\begin{equation}\label{Ikmn*}
    I^{k*}_{\mu\nu}:=(\nu,s+k,s+k-1,\cdots,\widehat \mu,\cdots, s-p+k+1)\,.
\end{equation}

It is clear that the coordinate $b_{i_1j_1}$ vanishes on $\mathfrak{D}$. Recalling the proof of Lemma 3.11 in \cite{F}, by (\ref{bst}) we can compute the defining equation $\rho_m$ for $\check B^{-j}_m$ in the local coordinate chart $\left(A^{\tau},(J_{j-1}^{\tau})^{-1}\right)$ as follows.
\begin{equation}
\begin{split}
&\rho_m=\frac{P_{I_m}\left(\Gamma_{j-1}^{\tau}\left(\widetilde X,\widetilde Y,\overrightarrow B^1,\cdots,\overrightarrow B^p\right)\right)}{\prod\limits_{t=1}^{m-j+1}b_{i_tj_t}^{m-j+2-t}}\,,\,\,\,j\leq m\leq r\,;\\
&\rho_{j-1}\equiv 1\,;\\
&\rho_m=\frac{P_{I_m}\left(\Gamma_{j-1}^{\tau}\left(\widetilde X,\widetilde Y,\overrightarrow B^1,\cdots,\overrightarrow B^p\right)\right)}{\prod\limits_{t=p-j+2}^{p-m}a_{i_{t}j_t}^{p-m-t+1}}\,,\,\,\,0\leq m\leq j-2\,.\\    
\end{split}
\end{equation}
In particular, we have that
\begin{equation}\label{nv}
    \rho_m(\mathfrak a)\neq 0\,,\,\,\,\,0\leq m\leq r\,.
\end{equation}
Therefore, by Claims II in the proof of Lemma 3.11 in \cite{F}, we have that
\smallskip

{\noindent\bf Property I.} For each $I\in\mathbb I_{s,p,n}^m$, $0\leq m\leq p$, the following rational function is nonzero at the point $\mathfrak a$,
\begin{equation}\label{ghol2}
\frac{P_{I}\left(\Gamma_{j-1}^{\tau}\left(\widetilde X,\widetilde Y,\overrightarrow B^1,\cdots,\overrightarrow B^p\right)\right)}{P_{I_m}\left(\Gamma_{j-1}^{\tau}\left(\widetilde X,\widetilde Y,\overrightarrow B^1,\cdots,\overrightarrow B^p\right)\right)}.
\end{equation} 

Denote by $\tau_0$ the index
\begin{equation}
\left(\begin{matrix}
j&j+1&\cdots&p&j-1&j-2&\cdots&1\\
s+j&s+j+1&\cdots&s+p&s-p+j-1&s-p+j-2&\cdots& s-p+1\\
\end{matrix}\right)\,.
\end{equation}
Denote the local coordinates of the holomorphic coordinate chart $\left(A^{\tau_0},(J_{j-1}^{\tau_0})^{-1}\right)$ as follows (see Appendix \ref{section:nc} as well). 
\begin{equation}
\widetilde X_0:=\left(\begin{matrix}
x_{0,1(s+j)}&\cdots &x_{0,1n}\\
\vdots&\ddots&\vdots\\
x_{0,l(s+j)}&\cdots &x_{0,ln}\\
\end{matrix}\right)\,\,\,{\rm and}\,\,\,\widetilde Y:=\left(\begin{matrix}
y_{0,j1}&\cdots& y_{0,j(s-p+j-1)}\\
\vdots&\ddots&\vdots\\
y_{0,p1}&\cdots& y_{0,p(s-p+j-1)}\\
\end{matrix}  \right)  \,;\,
\end{equation}
\begin{equation}
\begin{split}
&\overrightarrow B^{1}_0:=\left(b_{0,j(s+j)},\xi^{(1)}_{0,j(s+j+1)},\xi^{(1)}_{0,j(s+j+2)},\cdots,\xi^{(1)}_{0,jn},\xi^{(1)}_{0,(j+1)(s+j)},\xi^{(1)}_{0,(j+2)(s+j)},\cdots,\xi^{(1)}_{0,p(s+j)}\right)\,,\\
&\overrightarrow B^{2}_0:=\left(b_{0,(j+1)(s+j+1)},\xi^{(2)}_{0,(j+1)(s+j+2)},\xi^{(2)}_{0,(j+1)(s+j+3)},\cdots,\xi^{(2)}_{0,(j+1)n},\right.\,\\
&\,\,\,\,\,\,\,\,\,\,\,\,\,\,\,\,\,\,\,\,\,\,\,\,\,\,\,\,\,\,\left.\xi^{(2)}_{0,(j+2)(s+j+1)},\xi^{(2)}_{0,(j+3)(s+j+1)},\cdots,\xi^{(2)}_{0,p(s+j+1)}\right)\,,\\
&\,\,\,\,\,\,\,\,\,\,\,\,\,\,\,\,\,\,\,\,\,\,\,\,\,\,\,\,\vdots\\
&\overrightarrow B^{p-j+1}_0:=\left(b_{0,p(s+p)},\xi^{(p-j+1)}_{0,p(s+p+1)},\xi^{(p-j+1)}_{0,p(s+p+2)},\cdots,\xi^{(p-j+1)}_{0,pn}\right);\\
\end{split}
\end{equation}
\begin{equation}
\begin{split}
&\overrightarrow B^{p-j+2}_0:=\left(a_{0,(j-1)(s-p+j-1)},\xi^{(p-j+2)}_{0,(j-1)1},\xi^{(p-j+2)}_{0,(j-1)2},\cdots,\xi^{(p-j+2)}_{0,(j-1)(s-p+j-2)},\right.\,\,\,\,\,\,\,\,\,\,\,\,\,\,\,\,\,\,\,\,\,\,\,\,\,\,\,\,\,\,\,\,\,\,\,\,\,\,\,\,\,\,\,\,\,\,\,\,\,\,\,\,\,\,\,\,\,\\
&\,\,\,\,\,\,\,\,\,\,\,\,\,\,\,\,\,\,\,\,\,\,\,\,\,\,\,\,\,\,\,\,\,\,\,\,\,\,\,\,\,\left.\xi^{(p-j+2)}_{0,1(s-p+j-1)},\xi^{(p-j+2)}_{0,2(s-p+j-1)},\cdots,\xi^{(p-j+2)}_{0,(j-2)(s-p+j-1)}\right)\,,\,\,\,\,\,\,\,\,\,\,\,\,\,\\
&\overrightarrow B^{p-j+3}_0:=\left(a_{0,(j-2)(s-p+j-2)},\xi^{(p-j+3)}_{0,(j-2)1},\xi^{(p-j+3)}_{0,(j-2)2},\cdots,\xi^{(p-j+3)}_{0,(j-2)(s-p+j-2)},\right.\\
&\,\,\,\,\,\,\,\,\,\,\,\,\,\,\,\,\,\,\,\,\,\,\,\,\,\,\,\,\,\,\,\,\,\,\,\,\,\,\,\left.\xi^{(p-j+3)}_{0,1(s-p+j-2)},\xi^{(p-j+3)}_{0,2(s-p+j-2)},\cdots,\xi^{(p-j+3)}_{0,(j-3)(s-p+j-2)}\right)\,,\,\,\,\,\,\,\,\,\,\,\,\,\,\\
&\,\,\,\,\,\,\,\,\,\,\,\,\,\,\,\,\,\,\,\,\,\,\,\,\,\,\,\,\vdots\\
&\overrightarrow B^{p}_0:=\left(a_{0,1(s-p+1)},\xi^{(p)}_{0,11},\xi^{(p)}_{0,12},\cdots,\xi^{(p)}_{0,1(s-p)}\right)\,.\,\,\,\,\,\,\,\,\,\,\,\,\,\\
\end{split}
\end{equation}
Based on Claims III, III$^{\prime}$, III$^{\prime\prime}$, III$^{\prime\prime\prime}$,  and III$^{\prime\prime\prime\prime}$ in the proof of Lemma 3.11 in \cite{F}, we can make the change of coordinates between $\left(\widetilde X_0,\widetilde Y_0,\overrightarrow B^1_0,\cdots,\overrightarrow B^p_0\right)$ and $\left(\widetilde X,\widetilde Y,\overrightarrow B^1,\cdots,\overrightarrow B^p\right)$ as follows.
\begin{equation}
\small
\begin{split}
&b_{0,j(s+j)}=\frac{P_{I_{j}}\left(\Gamma_{j-1}^{\tau}\left(\widetilde X,\widetilde Y,\overrightarrow B^1,\cdots,\overrightarrow B^p\right)\right)}{P_{I_{j-1}}\left(\Gamma_{j-1}^{\tau}\left(\widetilde X,\widetilde Y,\overrightarrow B^1,\cdots,\overrightarrow B^p\right)\right)}\,;\\
&b_{0,m(s+m)}=\frac{P_{I_{m}}\left(\Gamma_{j-1}^{\tau}\left(\widetilde X,\widetilde Y,\overrightarrow B^1,\cdots,\overrightarrow B^p\right)\right)\cdot P_{I_{m-2}}\left(\Gamma_{j-1}^{\tau}\left(\widetilde X,\widetilde Y,\overrightarrow B^1,\cdots,\overrightarrow B^p\right)\right)}{\left(P_{I_{m-1}}\left(\Gamma_{j-1}^{\tau}\left(\widetilde X,\widetilde Y,\overrightarrow B^1,\cdots,\overrightarrow B^p\right)\right)\right)^2}\,,\,\,\,\,\,j+1\leq m\leq r=p\,;\\
&a_{0,(j-1)(s-p+j-1)}=\frac{P_{I_{j-2}}\left(\Gamma_{j-1}^{\tau}\left(\widetilde X,\widetilde Y,\overrightarrow B^1,\cdots,\overrightarrow B^p\right)\right)}{P_{I_{j-1}}\left(\Gamma_{j-1}^{\tau}\left(\widetilde X,\widetilde Y,\overrightarrow B^1,\cdots,\overrightarrow B^p\right)\right)}\,;\\
&a_{0,l(s-p+l)}=\frac{P_{I_{l+1}}\left(\Gamma_{j-1}^{\tau}\left(\widetilde X,\widetilde Y,\overrightarrow B^1,\cdots,\overrightarrow B^p\right)\right)\cdot P_{I_{l-1}}\left(\Gamma_{j-1}^{\tau}\left(\widetilde X,\widetilde Y,\overrightarrow B^1,\cdots,\overrightarrow B^p\right)\right)}{\left(P_{I_{l}}\left(\Gamma_{j-1}^{\tau}\left(\widetilde X,\widetilde Y,\overrightarrow B^1,\cdots,\overrightarrow B^p\right)\right)\right)^2}\,,\,\,\,\,\,1\leq l\leq j-2\,.\\
\end{split}
\end{equation}
For $1\leq m\leq p-j$,  $s-p+j+m\leq \mu\leq s$, and $\nu=s-p+k$, 
\begin{equation}
\xi^{(m)}_{0,(\mu-s+p)(s+j-1+m)}=\frac{P_{I_{\mu\nu}^{(j-1+m)}}\left(\Gamma_{j-1}^{\tau}\left(\widetilde X, \widetilde Y,\cdots,\overrightarrow B^p\right)\right)}{(-1)^{(j-1+m)(p-j+1-m)+\mu-s+p-j+1-m}\cdot P_{I_m}\left(\Gamma_{j-1}^{\tau}\left(\widetilde X,\widetilde Y,\overrightarrow B^1,\cdots,\overrightarrow B^p\right)\right)}\,;
\end{equation}
for $1\leq m\leq p-j+1$, $\mu=s+j-1+m$, and $s+j+m\leq \nu\leq n$, 
\begin{equation}\label{thol4}
\begin{split}
\xi^{(m)}_{0,(j-1+m)\nu}=\frac{P_{I_{\mu\nu}^{(j-1+m)*}}\left(\Gamma_{j-1}^{\tau}\left(\widetilde X, \widetilde Y,\cdots,\overrightarrow B^p\right)\right)}{(-1)^{(j-1+m)(p-j+1-m)}\cdot P_{I_m}\left(\Gamma_{j-1}^{\tau}\left(\widetilde X,\widetilde Y,\overrightarrow B^1,\cdots,\overrightarrow B^p\right)\right)}\,.\\
\end{split}
\end{equation}
For $p-j+2\leq m\leq p-1$, $s+1\leq \mu\leq s+p-m$, and $\nu=s+p-m+1$, 
\begin{equation}
\xi^{(m)}_{0,(\mu-s)(s-m+1)}=\frac{P_{I_{\mu\nu}^{(p-m)*}}\left(\Gamma_{j-1}^{\tau}\left(\widetilde X, \widetilde Y, \cdots,\overrightarrow B^p\right)\right)}{(-1)^{m(p-m)+s+p-m+1-\mu}\cdot P_{I_m}\left(\Gamma_{j-1}^{\tau}\left(\widetilde X,\widetilde Y,\overrightarrow B^1,\cdots,\overrightarrow B^p\right)\right)}\,;
\end{equation}
for $p-j+2\leq m\leq p$, $\mu=s-m+1$,  and $1\leq \nu\leq s-m$, 
\begin{equation}
\xi^{(m)}_{0,(p-m+1)\nu}=\frac{P_{I^{(p-m)}_{\mu\nu}}\left(\Gamma_{j-1}^{\tau}\left(\widetilde X, \widetilde Y,\cdots,\overrightarrow B^p\right)\right)}{(-1)^{m(p-m)} \cdot P_{I_m}\left(\Gamma_{j-1}^{\tau}\left(\widetilde X,\widetilde Y,\overrightarrow B^1,\cdots,\overrightarrow B^p\right)\right)}\,.
\end{equation}
For  $s-p+j\leq \mu\leq s$ and $1\leq \nu\leq s-p+j-1$, 
\begin{equation}
y_{0,(\mu-s+p)\nu}=(-1)^{k(p-k)+\mu-s+p-l}\cdot P_{I_{\mu\nu}^{(j-1)}}\left(\Gamma_{j-1}^{\tau}\left(\widetilde X, \widetilde Y,\cdots,\overrightarrow B^p\right)\right)\,;
\end{equation}
for  $s+1\leq \mu\leq s+j-1$, and $s+j\leq \nu\leq n$, 
\begin{equation}
x_{0,(\mu-s)\nu}=(-1)^{(j-1)(p-j+1)+\mu-s-j+1}\cdot P_{I_{\mu\nu}^{(j-1)*}}\left(\Gamma_{j-1}^{\tau}\left(\widetilde X, \widetilde Y,\cdots,\overrightarrow B^p\right)\right).
\end{equation}

Then by (\ref{nv}), it is easy to verify that the point $\mathfrak a$ has well-defined local coordinates in the holomorphic coordinate chart $\left(A^{\tau_0},(J_{j-1}^{\tau_0})^{-1}\right)$. 

We complete the proof of Claim. \,\,\,\,$\endpf$
\smallskip

According to Claim, it is easy to verify that the action of the Borel group $B$ on $\mathfrak a$ has a dense orbit in $D^-_j$. This is a contradiction. Hence, $\mathfrak D$ is an irreducible component of one of the divisors in $\{\check B^{- j}_{0}, \check B^{- j}_{1}, \cdots, \check B^{- j}_{r}\}$.

By a similar argument, we can prove (\ref{dbi2}) and (\ref{dbi3}) when $p=n-s<s$ and $p=n-s=s$, respectively.

We complete the proof of Lemma \ref{dgb}.\,\,\,\,$\endpf$ 
\begin{remark}
$\check B^{-1}_m$ (resp. $\check B^{+1}_m$) is irreducible if it is a proper subvariety of $D^{-1}$ (resp. $D^{+1}$). When $2\leq j\leq r$, $\check B^{-j}_m$ (resp. $\check B^{+j}_m$) consists of two components.
\end{remark}
\begin{proposition}\label{ccone}
The interior points of the cone in $A_{n-1}(D^{-}_j)\otimes_{\mathbb Z}\mathbb Q$ (resp. $A_{n-1}(D^{+}_j)\otimes_{\mathbb Z}\mathbb Q$) generated by  
\begin{equation}
\small
\begin{split}
&\left \{D^{-j}_{-1}, D^{-j}_{-2}, \cdots,\widehat{D^{-j}_{-j}}, \cdots, D^{-j}_{-r}, D^{-j}_{+(r+2-j)}, D^{-j}_{+(r+3-j)}, \cdots, D^{-j}_{+r},\check B^{-j}_{0}, \check B^{-j}_{1}, \cdots, \check B^{-j}_{r}\right\}\\
&\left(resp.\, \left \{D^{+j}_{+1}, D^{+j}_{+2}, \cdots,\widehat{D^{+j}_{+j}}, \cdots, D^{+j}_{+r}, D^{+j}_{-(r+2-j)}, D^{+j}_{-(r+3-j)}, \cdots, D^{+j}_{-r},\check B^{+j}_{0}, \check B^{+j}_{1}, \cdots, \check B^{+j}_{r}\right\}\right)\\
\end{split}
\end{equation}
are the interior points of the cone of effective divisors of $D^{-}_j$ (resp. $D^{+}_j$ ).
\end{proposition}

{\noindent\bf Proof of Proposition \ref{ccone}.}
Without loss of generality, we only consider the case $D^-_j$ in the following. Take the irreducible decomposition of $\check B^{-j}_m$ as follows.
\begin{equation}
\check B^{-j}_m=\sum_{\alpha=1}^{n^{-j}_m}\beta^{-j}_{m,\alpha}\cdot\Gamma^{-j}_{m,\alpha}, \,\,0\leq m\leq r,
\end{equation}
where $n^{-j}_m$ and $\beta^{-j}_{m,\alpha}$ are positive integers if 
$\check B^{-j}_m$ is nonempty. By Lemma \ref{cone}, the cone of effective divisors of $D^{-}_j$ is generated by
\begin{equation}
\begin{split}
 &\left\{D^{-j}_{-1}, D^{-j}_{-2}, \cdots,\widehat{D^{-j}_{-j}}, \cdots, D^{-j}_{-r}, D^{-j}_{+(r+2-j)}, D^{-j}_{+(r+3-j)}, \cdots, D^{-j}_{+r},\right.\\
 &\,\,\,\,\,\,\,\left.\Gamma^{-j}_{0,1},\cdots, \Gamma^{-j}_{0,n^{-j}_0},\Gamma^{-j}_{1,1},\cdots, \Gamma^{-j}_{1,n^{-j}_1},\cdots,\Gamma^{-j}_{r,1},\cdots, \Gamma^{-j}_{r,n^{-j}_r},\right\}.
\end{split}
\end{equation}
Then, Proposition \ref{ccone} follows.\,\,\,$\endpf$

\subsection{Bigness of the restriction of the anticanonical bundle}

Recall the following result.

\begin{lemma}[Lemma 6.7 in \cite {F}]\label{wk}
When  $p<n-s\leq s$,
\begin{equation}\label{tk1}
    K_{\mathcal T_{s,p,n}}=-(s-p+1)\cdot B_0 -2\sum_{j=1}^{p-1}B_j -(n-s-p+1)\cdot B_p -\sum_{i=1}^pD^-_i-\sum_{i=1}^pD^+_i\,;
\end{equation}
when $n-s=p<s$,
\begin{equation}\label{tk2}
    K_{\mathcal T_{s,p,n}}=-(s-p+1)\cdot B_0 -2\sum_{j=1}^{p-1}B_j -\sum_{i=1}^pD^-_i-\sum_{i=1}^pD^+_i\,;
\end{equation}
when $n-s<p<s$ ($r=n-s$),
\begin{equation}
    K_{\mathcal T_{s,p,n}}=-(s-p+1)\cdot B_0-2\sum_{j=1}^{r-1}B_j-(p-r+1)\cdot B_r-\sum_{i=1}^rD^-_i-\sum_{i=1}^rD^+_i\,;
\end{equation}
when $n-s=p=s$,
\begin{equation}
    K_{\mathcal T_{s,p,n}}= -2\sum_{j=1}^{p-1}B_j  -\sum_{i=1}^pD^-_i-\sum_{i=1}^pD^+_i\,.
\end{equation}
\end{lemma}

Computation yields that
\begin{lemma}\label{rwk}Let $1\leq j\leq r$. When  $r=p<n-s\leq s$, we have that
\begin{equation}
\begin{split}
&-K_{\mathcal T_{s,p,n}}|_{D^{- }_j}=(s-p+1)\cdot \check B^{-j}_0 +2\sum_{m=1}^{r-1}\check B^{-j}_m +(n-s-p+1)\cdot \check B^{-j}_r+\check B^{-j}_j-\check B^{-j}_{j+1}\\
&\,\,\,\,\,\,\,\,\,\,\,\,\,\,\,\,\,\,\,\,\,\,\,\,\,\,\,\,\,\,\,\,\,\,\,\,\,\,\,\,+\sum_{i=j+2}^rD^{-j}_{-i}+\sum_{i=r+2-j}^rD^{-j}_{+i}\,,\\
\end{split}
\end{equation}
\begin{equation}\label{d-j}
\sum_{i=1}^{j-1}D^{-j}_{-i}=\check B^{-j}_{j-2}-\check B^{-j}_{j-1}+D^{-j}_{+(r+2-j)}\,,\,\,\,\,\,\,\,2\leq j\leq r,\,\,\,\,\,\,\,\,\,\,\,\,\,\,\,\,\,\,\,\,\,\,\,\,\,\,\,\,\,\,\,\,\,\,\,\,\,\,\,\,\,\,\,\,\,\,\,\,\,\,\,\,\,\,\,\,\,\,\,\,\,\,\,\,\,\,\,\,\,
\end{equation}
\begin{equation}\label{d-j1}
D^{-j}_{-(j+1)} =-\check B^{-j}_{j-1}+2\check B^{-j}_{j}-\check B^{-j}_{j+1}\,,\,\,\,\,\,\,\,1\leq j\leq r-1.\,\,\,\,\,\,\,\,\,\,\,\,\,\,\,\,\,\,\,\,\,\,\,\,\,\,\,\,\,\,\,\,\,\,\,\,\,\,\,\,\,\,\,\,\,\,\,\,\,\,\,\,\,\,\,\,\,\,\,\,\,\,\,\,\,
\end{equation}
And
\begin{equation}
\begin{split}
-K_{\mathcal T_{s,p,n}}|_{D^{+ }_j}=&(s-p+1)\cdot \check B^{+j}_0 +2\sum_{m=1}^{r-1}\check B^{+j}_m +(n-s-p+1)\cdot \check B^{+j}_r+\check B^{+j}_{r-j}-\check B^{+j}_{r-j-1}\,,\\
&\,\,\,\,+\sum_{i=j+2}^rD^{+j}_{+i}+\sum_{i=r+2-j}^rD^{+j}_{-i}\\
\end{split}
\end{equation}
\begin{equation}\label{+d-j}
\sum_{i=1}^{j-1}D^{+j}_{+i}=\check B^{+j}_{r+2-j}-\check B^{+j}_{r+1-j}+D^{+j}_{-(r+2-j)}\,,\,\,\,\,\,\,2\leq j\leq r,\,\,\,\,\,\,\,\,\,\,\,\,\,\,\,\,\,\,\,\,\,\,\,\,\,\,\,\,\,\,\,\,\,\,\,\,\,\,\,\,\,\,\,\,\,\,\,\,\,\,\,\,\,\,\,\,\,\,\,\,\,\,\,\,
\end{equation}
\begin{equation}\label{+d-j1}
D^{+j}_{+(j+1)} =-\check B^{+j}_{r+1-j}+2\check B^{+j}_{r-j}-\check B^{+j}_{r-1-j}\,,\,\,\,\,\,\,\,1\leq j\leq r-1.\,\,\,\,\,\,\,\,\,\,\,\,\,\,\,\,\,\,\,\,\,\,\,\,\,\,\,\,\,\,\,\,\,\,\,\,\,\,\,\,\,\,\,\,\,\,\,\,\,\,\,\,\,\,\,\,\,
\end{equation}
Notice that here we use the convention that $\check B^{-r}_{r+1}$ and $\check B^{+r}_{-1}$ are trivial.
\end{lemma}

{\noindent\bf Proof of Lemma \ref{rwk}.}
Since $D^+_{1},D^+_2,\cdots,D^+_{r+1-j}$ have empty intersection with $D^-_j$, by Lemma \ref{bb} we can conclude the following formulas for line bundles.
\begin{equation}\label{dbst}
\begin{split}
\check B^{-j}_0&=(R_{s,p,n})^*\left(\mathcal O_{G(p,n)}(1)\right)|_{D^-_j}-\sum_{i=r+2-j}^{r}(r+1-i)\cdot D^{-j}_{+i}\,,\\
\check B^{-j}_1&=(R_{s,p,n})^*\left(\mathcal O_{G(p,n)}(1)\right)|_{D^-_j}-\sum_{i=r+2-j}^{r-1}(r-i)\cdot D^{-j}_{+i}-D^{-j}_{-1}\,,\\
\check B^{-j}_2&=(R_{s,p,n})^*\left(\mathcal O_{G(p,n)}(1)\right)|_{D^-_j}-\sum_{i=r+2-j}^{r-2}(r-1-i)\cdot D^{-j}_{+i}-2D^{-j}_{-1}-D^{-j}_{-2}\,,\\
&\,\,\,\vdots\\
\end{split}
\end{equation}
\begin{equation}\label{dm-2}
\begin{split}
\check B^{-j}_{j-2}&=(R_{s,p,n})^*\left(\mathcal O_{G(p,n)}(1)\right)|_{D^-_j}-D^{-j}_{+(r+2-j)}-(j-2)D^{-j}_{-1}-(j-3)D^{-j}_{-2}-\cdots-D^{-j}_{-(j-2)}\,,\\
\end{split}
\end{equation}
\begin{equation}\label{dm-1}
\begin{split}
\check B^{-j}_{j-1}&=(R_{s,p,n})^*\left(\mathcal O_{G(p,n)}(1)\right)|_{D^-_j}-(j-1)D^{-j}_{-1}-(j-2)D^{-j}_{-2}-\cdots-D^{-j}_{-(j-1)}\,,\\
\end{split}
\end{equation}
\begin{equation}\label{dm}
\begin{split}
\check B^{-j}_j&=(R_{s,p,n})^*\left(\mathcal O_{G(p,n)}(1)\right)|_{D^-_j}-jD^{-j}_{-1}-(j-1)D^{-j}_{-2}-\cdots-2D^{-j}_{-(j-1)}-D^{-j}_{-j}\,,\\
\end{split}
\end{equation}
\begin{equation}\label{dm+1}
\begin{split}
\check B^{-j}_{j+1}&=(R_{s,p,n})^*\left(\mathcal O_{G(p,n)}(1)\right)|_{D^-_j}-(j+1)D^{-j}_{-1}-jD^{-j}_{-2}-\cdots-2D^{-j}_{-j}-D^{-j}_{-(j+1)}\,,\\
\end{split}
\end{equation}
\begin{equation}
\begin{split}
&\,\,\,\vdots\\
\check B^{-j}_r&=(R_{s,p,n})^*\left(\mathcal O_{G(p,n)}(1)\right)|_{D^-_j}-rD^{-j}_{-1}-\cdots-(r+1-j)D^{-j}_{-j}-\cdots-D^{-j}_{-r}\,.\\
\end{split}
\end{equation} 

Then (\ref{d-j}) and (\ref{d-j1}) follow from (\ref{dm-2}), (\ref{dm-1}), (\ref{dm}), and (\ref{dm+1}) directly.

Subtracting (\ref{dm+1}) from (\ref{dm}), we derive that
\begin{equation}\label{djj}
D^{-j}_{-1}+D^{-j}_{-2}+\cdots+D^{-j}_{-(j-1)}+D^{-j}_{-j}+D^{-j}_{-(j+1)}=\check B^{-j}_j-\check B^{-j}_{j+1}.\\
\end{equation}
Restricting (\ref{tk1}) to $D^-_j$ and plugging in (\ref{djj}), we have that
\begin{equation}
\begin{split}
 -K_{\mathcal T_{s,p,n}}|_{D^{- }_j}=&(s-p+1)\cdot \check B^{-j}_0 +2\sum_{m=1}^{r-1}\check B^{-j}_m +(n-s-p+1)\cdot \check B^{-j}_r\\
&\,\,\,\,\,+\sum_{i=1}^{j+1}D^{-j}_{-i}+
\sum_{i=j+2}^rD^{-j}_{-i}+\sum_{i=r+2-j}^rD^{-j}_{+i}\,\\
=&(s-p+1)\cdot \check B^{-j}_0 +2\sum_{m=1}^{r-1}\check B^{-j}_m +(n-s-p+1)\cdot \check B^{-j}_r+\check B^{-j}_j-\check B^{-j}_{j+1}\\
&\,\,\,\,\,+\sum_{i=j+2}^{r}D^{-j}_{-i}+\sum_{i=r+2-j}^rD^{-j}_{+i}\,.\\
\end{split}
\end{equation}

Similarly, by Lemma \ref{bb} and the fact that $D^-_{1},D^-_2,\cdots,D^-_{r+1-j}$ have empty intersection with $D^+_j$, we can conclude that
\begin{equation}\label{+dbst}
\begin{split}
\check B^{+j}_0&=(R_{s,p,n})^*\left(\mathcal O_{G(p,n)}(1)\right)|_{D^+_j}-\sum_{i=1}^{r}(r+1-i)\cdot D^{+j}_{+i}\,,\\
\check B^{+j}_1&=(R_{s,p,n})^*\left(\mathcal O_{G(p,n)}(1)\right)|_{D^+_j}-\sum_{i=1}^{r-1}(r-i)\cdot D^{+j}_{+i}\,,\\
&\,\,\,\vdots\\
\end{split}
\end{equation}
\begin{equation}\label{+dm-2}
\begin{split}
\check B^{+j}_{r-1-j}&=(R_{s,p,n})^*\left(\mathcal O_{G(p,n)}(1)\right)|_{D^+_j}-\sum_{i=1}^{j+1}(j+2-i)\cdot D^{+j}_{+i}\,,\\
\end{split}
\end{equation}
\begin{equation}\label{+dm-1}
\begin{split}
\check B^{+j}_{r-j}&=(R_{s,p,n})^*\left(\mathcal O_{G(p,n)}(1)\right)|_{D^+_j}-\sum_{i=1}^{j}(j+1-i)\cdot D^{+j}_{+i}\,,\\
\end{split}
\end{equation}
\begin{equation}\label{+dm}
\begin{split}
\check B^{+j}_{r+1-j}&=(R_{s,p,n})^*\left(\mathcal O_{G(p,n)}(1)\right)|_{D^+_j}-\sum_{i=1}^{j-1}(j-i)\cdot D^{+j}_{+i}\,,\\
\end{split}
\end{equation}
\begin{equation}\label{+dm+1}
\begin{split}
\check B^{+j}_{r+2-j}&=(R_{s,p,n})^*\left(\mathcal O_{G(p,n)}(1)\right)|_{D^+_j}-\sum_{i=1}^{j-2}(j-1-i)\cdot D^{+j}_{+i}-D^{+j}_{-(r+2-j)}\,,\\
\end{split}
\end{equation}
\begin{equation}
\begin{split}
&\,\,\,\vdots\\
\check B^{+j}_r&=(R_{s,p,n})^*\left(\mathcal O_{G(p,n)}(1)\right)|_{D^+_j}-\sum_{i=r+2-j}^{r}(r+1-i)\cdot D^{+j}_{-i}\,.\\
\end{split}
\end{equation} 

Then (\ref{+d-j}) and (\ref{+d-j1}) follow from (\ref{+dm-2}), (\ref{+dm-1}), (\ref{+dm}), and (\ref{+dm+1}) directly.

Subtracting (\ref{+dm-2}) from (\ref{+dm-1}), we derive that
\begin{equation}\label{+djj}
D^{+j}_{+1}+D^{+j}_{+2}+\cdots+D^{+j}_{+(j-1)}+D^{+j}_{+j}+D^{+j}_{+(j+1)}=\check B^{+j}_{r-j}-\check B^{+j}_{r-1-j}.\\
\end{equation}
Restricting (\ref{tk1}) to $D^+_j$ and plugging in (\ref{+djj}), we have that
\begin{equation}
\begin{split}
 -K_{\mathcal T_{s,p,n}}|_{D^{+}_j}=&(s-p+1)\cdot \check B^{+j}_0 +2\sum_{m=1}^{r-1}\check B^{+j}_m +(n-s-p+1)\cdot \check B^{+j}_r\\
&\,\,\,\,\,+\sum_{i=1}^{j+1}D^{+j}_{+i}+
\sum_{i=j+2}^rD^{+j}_{+i}+\sum_{i=r+2-j}^rD^{+j}_{-i}\,,\\
=&(s-p+1)\cdot \check B^{+j}_0 +2\sum_{m=1}^{r-1}\check B^{+j}_m +(n-s-p+1)\cdot \check B^{+j}_r+\check B^{+j}_{r-j}-\check B^{+j}_{r-1-j}\\
&\,\,\,\,\,+\sum_{i=j+2}^{r}D^{+j}_{+i}+\sum_{i=r+2-j}^rD^{+j}_{-i}\,.\\
\end{split}
\end{equation}

We complete the proof of Lemma \ref{rwk}.\,\,\,\,$\endpf$ 
\medskip

\begin{lemma}\label{rwk1}Let $1\leq j\leq r$. When  $r=n-s<p<s$, we have that
\begin{equation}
\begin{split}
&-K_{\mathcal T_{s,p,n}}|_{D^{- }_j}=(s-p+1)\cdot \check B^{-j}_0 +2\sum_{m=1}^{r-1}\check B^{-j}_m +(p-r+1)\cdot \check B^{-j}_r+\check B^{-j}_j-\check B^{-j}_{j+1}\\
&\,\,\,\,\,\,\,\,\,\,\,\,\,\,\,\,\,\,\,\,\,\,\,\,\,\,\,\,\,\,\,\,\,\,\,\,\,\,\,\,+\sum_{i=j+2}^rD^{-j}_{-i}+\sum_{i=r+2-j}^rD^{-j}_{+i}\,,\\
\end{split}
\end{equation}
\begin{equation}\label{d-j-1}
\sum_{i=1}^{j-1}D^{-j}_{-i}=\check B^{-j}_{j-2}-\check B^{-j}_{j-1}+D^{-j}_{+(r+2-j)}\,,\,\,\,\,\,\,\,2\leq j\leq r,\,\,\,\,\,\,\,\,\,\,\,\,\,\,\,\,\,\,\,\,\,\,\,\,\,\,\,\,\,\,\,\,\,\,\,\,\,\,\,\,\,\,\,\,\,\,\,\,\,\,\,\,\,\,\,\,\,\,\,\,\,\,\,\,\,\,\,\,\,
\end{equation}
\begin{equation}\label{d-j1-1}
D^{-j}_{-(j+1)} =-\check B^{-j}_{j-1}+2\check B^{-j}_{j}-\check B^{-j}_{j+1}\,,\,\,\,\,\,\,\,1\leq j\leq r-1.\,\,\,\,\,\,\,\,\,\,\,\,\,\,\,\,\,\,\,\,\,\,\,\,\,\,\,\,\,\,\,\,\,\,\,\,\,\,\,\,\,\,\,\,\,\,\,\,\,\,\,\,\,\,\,\,\,\,\,\,\,\,\,\,\,
\end{equation}
And
\begin{equation}
\begin{split}
-K_{\mathcal T_{s,p,n}}|_{D^{+ }_j}=&(s-p+1)\cdot \check B^{+j}_0 +2\sum_{m=1}^{r-1}\check B^{+j}_m +(p-r+1)\cdot \check B^{+j}_r+\check B^{+j}_{r-j}-\check B^{+j}_{r-j-1}\\
&\,\,\,\,+\sum_{i=j+2}^rD^{+j}_{+i}+\sum_{i=r+2-j}^rD^{+j}_{-i}\,,\\
\end{split}
\end{equation}
\begin{equation}\label{+d-j-1}
\sum_{i=1}^{j-1}D^{+j}_{+i}=\check B^{+j}_{r+2-j}-\check B^{+j}_{r+1-j}+D^{+j}_{-(r+2-j)}\,,\,\,\,\,\,\,2\leq j\leq r,\,\,\,\,\,\,\,\,\,\,\,\,\,\,\,\,\,\,\,\,\,\,\,\,\,\,\,\,\,\,\,\,\,\,\,\,\,\,\,\,\,\,\,\,\,\,\,\,\,\,\,\,\,\,\,\,\,\,\,\,\,\,\,\,
\end{equation}
\begin{equation}\label{+d-j1-1}
D^{+j}_{+(j+1)} =-\check B^{+j}_{r+1-j}+2\check B^{+j}_{r-j}-\check B^{+j}_{r-1-j}\,,\,\,\,\,\,\,\,1\leq j\leq r-1.\,\,\,\,\,\,\,\,\,\,\,\,\,\,\,\,\,\,\,\,\,\,\,\,\,\,\,\,\,\,\,\,\,\,\,\,\,\,\,\,\,\,\,\,\,\,\,\,\,\,\,\,\,\,\,\,\,
\end{equation}
Notice that by convention $\check B^{-r}_{r+1}$ and $\check B^{+r}_{-1}$ are trivial.
\end{lemma}
{\noindent\bf Proof of Lemma \ref{rwk1}.}
The proof is exactly the same as that of Lemma \ref{rwk1} by setting $r=n-s$ instead of $r=p$. For simplicity, we omit it here.\,\,\,\,$\endpf$

\begin{lemma}\label{rwk2} Let  $r=n-s=p<s$. we have the following identities. When $1\leq j\leq r-2$, 
\begin{equation}
\begin{split}
&-K_{\mathcal T_{s,p,n}}|_{D^{- }_j}=(s-p+1)\cdot \check B^{-j}_0 +2\sum_{m=1}^{r-1}\check B^{-j}_m+\check B^{-j}_j-\check B^{-j}_{j+1}+\sum_{i=j+2}^rD^{-j}_{-i}+\sum_{i=r+2-j}^rD^{-j}_{+i}\,;\\
\end{split}
\end{equation}
when $j=r-1$ or $r$,
\begin{equation}
\begin{split}
&-K_{\mathcal T_{s,p,n}}|_{D^{- }_j}=(s-p+1)\cdot \check B^{-j}_0 +2\sum_{m=1}^{r-1}\check B^{-j}_m+\check B^{-j}_{r-1}+\sum_{i=j+2}^rD^{-j}_{-i}+\sum_{i=r+2-j}^rD^{-j}_{+i}\,.\\
\end{split}
\end{equation}
\begin{equation}\label{1rd-j-1}
\sum_{i=1}^{j-1}D^{-j}_{-i}=\check B^{-j}_{j-2}-\check B^{-j}_{j-1}+D^{-j}_{+(r+2-j)}\,,\,\,\,\,\,\,\,2\leq j\leq r\,;\,\,\,\,\,\,\,\,\,\,\,\,\,\,\,\,\,\,\,\,\,\,\,\,\,\,\,\,\,\,\,\,\,\,\,\,\,\,\,\,\,\,\,\,\,\,\,\,\,\,\,\,\,\,\,\,\,\,\,\,\,\,\,\,\,\,\,\,
\end{equation}
\begin{equation}\label{1rd-j1-1}
D^{-j}_{-(j+1)} =\left\{\begin{matrix}
-\check B^{-j}_{j-1}+2\check B^{-j}_{j}-\check B^{-j}_{j+1}\,,&1\leq j\leq r-2\,\\
-\check B^{-j}_{j-1}+2\check B^{-j}_{j}\,,&j=r-1\,\\
\end{matrix}\right..
\,\,\,\,\,\,\,\,\,\,\,\,\,\,\,\,\,\,\,\,\,\,\,\,\,\,\,\,\,\,\,\,\,\,\,\,\,\,\,\,\,\,\,\,\,\,\,\,\,\,\,\,\,\,\,\,\,\,\,\,
\end{equation}
Similarly, when $1\leq j\leq r-1$,
\begin{equation}
\begin{split}
-K_{\mathcal T_{s,p,n}}|_{D^{+ }_j}=&(s-p+1)\cdot \check B^{+j}_0 +2\sum_{m=1}^{r-1}\check B^{+j}_m +\check B^{+j}_{r-j}-\check B^{+j}_{r-1-j}+\sum_{i=j+2}^rD^{+j}_{+i}+\sum_{i=r+2-j}^rD^{+j}_{-i}\,;\\
\end{split}
\end{equation}
when $j=r$,
\begin{equation}
\begin{split}
-K_{\mathcal T_{s,p,n}}|_{D^{+ }_j}=&(s-p+1)\cdot \check B^{+j}_0 +2\sum_{m=1}^{r-1}\check B^{+j}_m +\check B^{+j}_{1}-\check B^{+j}_{0}+\sum_{i=j+2}^rD^{+j}_{+i}+\sum_{i=r+2-j}^rD^{+j}_{-i}\,.\\
\end{split}
\end{equation}
\begin{equation}\label{r+d-j-1}
\sum_{i=1}^{j-1}D^{+j}_{+i}=\left\{\begin{matrix}\check B^{+j}_{r+2-j}-\check B^{+j}_{r+1-j}+D^{+j}_{-(r+2-j)}\,,&\,\,\,\,\,\,3\leq j\leq r\\
-\check B^{+j}_{r+1-j}+D^{+j}_{-(r+2-j)}\,,\,\,\,\,\,\,\,\,&\,\,\,\,\,\,j=2\\
\end{matrix}\,;\right.\,\,\,\,\,\,\,\,\,\,\,\,\,\,\,\,\,\,\,\,\,\,\,\,\,\,\,\,\,\,\,\,\,\,\,\,\,\,\,\,\,\,\,\,\,\,\,\,\,\,\,\,\,\,\,\,\,\,\,\,\,\,\,\,
\end{equation}
\begin{equation}\label{r+d-j1-1}
D^{+j}_{+(j+1)} =\left\{\begin{matrix}
-\check B^{+j}_{r+1-j}+2\check B^{+j}_{r-j}-\check B^{+j}_{r-1-j}\,,&\,\,\,\,\,\,\,2\leq j\leq r-1\\
-D^{+j}_{-r}+2\check B^{+j}_{r-j}-\check B^{+j}_{r-1-j}\,,&\,\,\,\,\,\,\,j=1\\
\end{matrix}.\right.\,\,\,\,\,\,\,\,\,\,\,\,\,\,\,\,\,\,\,\,\,\,\,\,\,\,\,\,\,\,\,\,\,\,\,\,\,\,\,\,\,\,\,\,\,\,\,\,\,\,\,\,\,\,\,\,\,
\end{equation}
Notice that by convention $\check B^{-r}_{r+1}$ and $\check B^{+r}_{-1}$ are  trivial line bundles.
\end{lemma}

{\noindent\bf Proof of Lemma \ref{rwk2}.} See Appendix \ref{dvanderl}. \,\,\,\,$\endpf$

\begin{lemma}\label{rwk3}
Let $r=n-s=p=s$. When $1\leq j\leq r-2$,
\begin{equation}
\begin{split}
-K_{\mathcal T_{s,p,n}}|_{D^{- }_j}=&2\sum_{m=1}^{r-1}\check B^{-j}_m+\check B^{-j}_j-\check B^{-j}_{j+1}+\sum_{i=j+2}^rD^{-j}_{-i}+\sum_{i=r+2-j}^rD^{-j}_{+i}\,;\\
\end{split}
\end{equation}
 when $j=r-1$ or $r$,
\begin{equation}
\begin{split}
-K_{\mathcal T_{s,p,n}}|_{D^{- }_j}=&2\sum_{m=1}^{r-1}\check B^{-j}_m+\check B^{-j}_{r-1}+\sum_{i=j+2}^rD^{-j}_{-i}+\sum_{i=r+2-j}^rD^{-j}_{+i}\,.\\
\end{split}
\end{equation}
\begin{equation}\label{sum1}
\sum_{i=1}^{j-1}D^{-j}_{-i}=\left\{\begin{matrix}\check B^{-j}_{j-2}-\check B^{-j}_{j-1}+D^{-j}_{+(r+2-j)}\,,&3\leq j\leq r\\
-\check B^{-j}_{j-1}+D^{-j}_{+(r+2-j)}\,,\,\,\,\,\,\,\,\,\,\,\,\,\,\,\,\,&j=2\\
\end{matrix}\right..
\end{equation}
\begin{equation}\label{sum2}
D^{-j}_{-(j+1)} =\left\{\begin{matrix}
-D^{-j}_{+r}+2\check B^{-j}_{j}-\check B^{-j}_{j+1}\,,&j=1\\
-\check B^{-j}_{j-1}+2\check B^{-j}_{j}-\check B^{-j}_{j+1}\,,&2\leq j\leq r-2\\
-\check B^{-j}_{j-1}+2\check B^{-j}_{j}\,,\,\,\,\,\,\,\,\,\,\,\,\,&j=r-1\\
\end{matrix}\right.
\end{equation}
Similarly, when $1\leq j\leq r-2$,
\begin{equation}
-K_{\mathcal T_{s,p,n}}|_{D^{+ }_j}=2\sum_{m=1}^{r-1}\check B^{+j}_m+\check B^{+j}_{r-j}-\check B^{+j}_{r-j-1}+\sum_{i=j+2}^rD^{+j}_{+i}+\sum_{i=r+2-j}^rD^{+j}_{-i};\\    
\end{equation}
when $j=r-1$ or $r$,
\begin{equation}
 -K_{\mathcal T_{s,p,n}}|_{D^{+ }_j}=2\sum_{m=1}^{r-1}\check B^{+j}_m +\check B^{+j}_{1}+\sum_{i=j+2}^{r}D^{+j}_{+i}+\sum_{i=r+2-j}^rD^{+j}_{-i}\,.\\
\end{equation}
\begin{equation}\label{pp1}
\sum_{i=1}^{j-1}D^{+j}_{+i}=\left\{\begin{matrix}\check B^{+j}_{r+2-j}-\check B^{+j}_{r+1-j}+D^{+j}_{-(r+2-j)}\,,&3\leq j\leq r\\
-\check B^{+j}_{r+1-j}+D^{-j}_{+(r+2-j)}\,,\,\,\,\,\,\,\,\,\,\,\,\,\,\,\,\,&j=2\\
\end{matrix}\right.;
\end{equation}
\begin{equation}\label{pp2}
D^{+j}_{+(j+1)} =\left\{\begin{matrix}
-\check B^{+j}_{r-j-1}+2\check B^{+j}_{r-j}-D^{+j}_{-r}\,,&j=1\\
-\check B^{+j}_{r-j-1}+2\check B^{+j}_{r-j}-\check B^{+j}_{r-j+1}\,,&2\leq j\leq r-2\\
2\check B^{+j}_{r-j}-\check B^{+j}_{r-j+1}\,,\,\,\,\,\,\,\,\,\,\,\,\,&j=r-1\\
\end{matrix}\right..
\end{equation}
\end{lemma}

{\noindent\bf Proof of Lemma \ref{rwk3}.} See Appendix \ref{dvanderld}.\,\,\,\,$\endpf$

\begin{proposition}\label{big}
When $r\geq 3$ and $1\leq i\leq r$, the restriction of the anticanonical bundle $-K_{\mathcal T_{s,p,n}}$ on $D^{\pm}_{i}$ is big and nef.
\end{proposition}
{\noindent\bf Proof of Proposition \ref{big}.}
By Theorem \ref{fano}, $-K_{\mathcal T_{s,p,n}}$ is nef. It is clear that the restriction $-K_{\mathcal T_{s,p,n}}|_{D^{\pm}_{i}}$ on $D^{\pm}_{i}$ is nef as well.

Since $D^{\pm}_i$ is projective, each interior point of the cone of effective divisors of $D^{\pm}_j$ represents a big divisor. Therefore, by Proposition \ref{ccone}, it suffices to show that $-K_{\mathcal T_{s,p,n}}|_{D^{\pm}_{i}}$ can be written as a positive combination of the following divisors.
\begin{equation}
   \left \{D^{\pm j}_{\pm 1}, D^{\pm j}_{\pm 2}, \cdots,\widehat{D^{\pm j}_{\pm j}}, \cdots, D^{\pm j}_{\pm r}, D^{\pm j}_{\mp (r+2-j)}, D^{\pm j}_{\mp (r+3-j)}, \cdots, D^{\pm j}_{\mp r},\check B^{\pm j}_{0}, \check B^{\pm j}_{1}, \cdots, \check B^{\pm j}_{r}\right\}
\end{equation}

Assume that $r=p<n-s\leq s$. By Lemma \ref{rwk}, we have that
\begin{equation}
\begin{split}
&-K_{\mathcal T_{s,p,n}}|_{D^{- }_j}=(s-p+1)\cdot \check B^{-j}_0 +2\sum_{m=1}^{r-1}\check B^{-j}_m +(n-s-p+1)\cdot \check B^{-j}_r+\check B^{-j}_j-\check B^{-j}_{j+1}\\
&\,\,\,\,\,\,\,\,\,\,\,\,\,\,\,\,\,\,\,\,\,\,\,\,\,\,\,\,\,\,\,\,\,\,\,\,\,\,\,\,+\sum_{i=j+2}^rD^{-j}_{-i}+\sum_{i=r+2-j}^rD^{-j}_{+i}\,,\\
&\,\,\,\,\,\,\,\,\,\,\,\,\,\,\,\,\,\,\,\,\,\,\,\,\,\,\,\,\,\,\,\,\,\,\,\,\,\,\,\,+\delta_1\left(\sum_{i=1}^{j-1}D^{-j}_{-i}-\left(\check B^{-j}_{j-2}-\check B^{-j}_{j-1}+D^{-j}_{+(r+2-j)}\right)\right)\\
&\,\,\,\,\,\,\,\,\,\,\,\,\,\,\,\,\,\,\,\,\,\,\,\,\,\,\,\,\,\,\,\,\,\,\,\,\,\,\,\,+\delta_2\left(D^{-j}_{-(j+1)}-\left(-\check B^{-j}_{j-1}+2\check B^{-j}_{j}-\check B^{-j}_{j+1}\right)\right)\,.\\
\end{split}  
\end{equation}
By choosing suitable $\delta_1,\delta_2$, we can conclude that $-K_{\mathcal T_{s,p,n}}|_{D^{- }_j}$ is in the interior of the cone of effective divisors of $D^{- }_j$ for $1\leq j\leq r$. Similarly,
\begin{equation}
\begin{split}
-K_{\mathcal T_{s,p,n}}|_{D^{+ }_j}=&(s-p+1)\cdot \check B^{+j}_0 +2\sum_{m=1}^{r-1}\check B^{+j}_m +(n-s-p+1)\cdot \check B^{+j}_r+\check B^{+j}_{r-j}-\check B^{+j}_{r-j-1}\,,\\
&\,\,\,\,+\sum_{i=j+2}^rD^{+j}_{+i}+\sum_{i=r+2-j}^rD^{+j}_{-i}\\
&\,\,\,\,+\delta_1\left(\sum_{i=1}^{j-1}D^{+j}_{+i}-\left(\check B^{+j}_{r+2-j}-\check B^{+j}_{r+1-j}+D^{+j}_{-(r+2-j)}\right)\right)\\
&\,\,\,\,+\delta_2\left(D^{+j}_{+(j+1)}-\left(\check B^{+j}_{r+1-j}+2\check B^{+j}_{r-j}-\check B^{+j}_{r-1-j}\right)\right)\,.
\end{split}
\end{equation}
By choosing suitable $\delta_1,\delta_2$, we can show that $-K_{\mathcal T_{s,p,n}}|_{D^{+}_j}$ is in the interior of the cone of effective divisors of $D^{- }_j$ for $1\leq j\leq r$. 

Similarly, we can prove that $-K_{\mathcal T_{s,p,n}}|_{D^{\pm}_j}$ is big for $1\leq j\leq r$ when $r=n-s<p$, $r=n-s=p<s$, or $r=n-s=p=s$.

We complete the proof of Proposition \ref{big}.
\subsection{Proof of Theorem \ref{bignef}.}


The proof is similar to that of Proposition 4.2 in \cite{BB}. When $r\leq 2$, Theorem \ref{bignef} follows from \cite{BB} directly, for $\mathcal T_{s,p,n}$ is regular by Lemma \ref{regu} and ample by Proposition \ref{fano}. In what follows, we will assume that $r\geq 3$. 

Recall the notion of the action sheaf of a spherical variety as follows.  Let $G$ be a connected reductive algebraic group over $\mathbb C$. Let $X$ be a spherical $G$-variety with boundary $\partial X$, that is, $\partial X=X\backslash\Omega$ where $\Omega$ is the open $G$-orbit of $X$. The action sheaf $S_X$ of $X$ is the subsheaf of $T_X$ made of vector fields tangent to $\partial X$. Combining Theorem 4.1 in \cite{Kn} and Proposition 2.5 in \cite{BB}, one can show that $H^i(X, S_X)=0$ for any complete regular variety and any $i>0$.  By Lemma \ref{regu}, we have that \begin{equation}\label{act}
H^i(\mathcal T_{s,p,n}, S_{\mathcal T_{s,p,n}})=0\,,\,\,\,\,i>0.
\end{equation}

By Proposition 2.3.2 in \cite{BB}, the following exact sequence of sheaves holds.
\begin{equation}\label{exa}
    0\rightarrow\mathcal{S}_{\mathcal T_{s,p,n}}\rightarrow \mathcal{T}_{\mathcal T_{s,p,n}}\rightarrow\bigoplus_{i=1}^r\mathcal O_{\mathcal T_{s,p,n}}\left(D_i^-\right)\otimes_{\mathcal O\left(\mathcal T_{s,p,n}\right)}\mathcal O_{D_i^-}\bigoplus_{i=1}^r\mathcal O_{\mathcal T_{s,p,n}}\left(D_i^+\right)\otimes_{\mathcal O\left(\mathcal T_{s,p,n}\right)}\mathcal O_{D_i^+}\rightarrow 0
\end{equation}
Since $D^{\pm}_i$ is smooth, $\mathcal O_{\mathcal T_{s,p,n}}\left(D_i^{\pm}\right)\otimes_{\mathcal O\left(\mathcal T_{s,p,n}\right)}\mathcal O_{D_i^{\pm}}\cong\mathcal N_{i}^{\pm}$ where $\mathcal N_{i}^{\pm}$ is the normal bundle of $D^{\pm}_i$ in $\mathcal T_{s,p,n}$. By the adjunction formula, $\mathcal N^{\pm}_i=-K_{\mathcal T_{s,p,n}}|_{D^{\pm}_i}+K_{D^{\pm}_i}$. By Proposition \ref{big}, $-K_{\mathcal T_{s,p,n}}|_{D^{\pm}_i}$ is big and nef. Then the Kawamata-Viehweg vanishing theorem yields that \begin{equation}\label{nor}
    H^j(D^{\pm}_i,\mathcal N^{\pm}_i)=0\,,\,\,j>0.
\end{equation}

Taking the long exact sequence of  (\ref{exa}), we can conclude by (\ref{act}) and (\ref{nor}) that
\begin{equation}
    H^j(\mathcal T_{s,p,n}, T_{\mathcal T_{s,p,n}})=\left(\bigoplus_{i=1}^rH^j(D^-_i, \mathcal N^-_{i})\right)\bigoplus\left(\bigoplus_{i=1}^rH^j(D^+_i, \mathcal N^{+}_{i})\right)=0\,,\,\,j>0\,.
\end{equation}

We complete the proof of Theorem \ref{bignef}.
\,\,\,\,$\endpf$
\medskip

\appendix

\section{Holomorphic Atlas}\label{section:nc}
Assume $n\leq 2s$, $p\leq \frac{n}{2}$, and $l\leq \min \{n-s,p\}$. Let $U_l$ be an affine open subset of $G(p,n)$ defined by
\begin{equation}\label{ul}
U_l:=\left\{\,\,\,\,\underbracedmatrixll{Z\\Y}{s-p+l\,\,\rm columns}
  \hspace{-.45in}\begin{matrix}
  &\hfill\tikzmark{a}\\
  &\hfill\tikzmark{b}  
  \end{matrix} \,\,\,\,\,
  \begin{matrix}
  0\\
I_{(p-l)\times(p-l)}\\
\end{matrix}\hspace{-.11in}
\begin{matrix}
  &\hfill\tikzmark{c}\\
  &\hfill\tikzmark{d}
  \end{matrix}\hspace{-.11in}\begin{matrix}
  &\hfill\tikzmark{g}\\
  &\hfill\tikzmark{h}
  \end{matrix}\,\,\,\,
\begin{matrix}
I_{l\times l}\\
0\\
\end{matrix}\hspace{-.11in}
\begin{matrix}
  &\hfill\tikzmark{e}\\
  &\hfill\tikzmark{f}\end{matrix}\hspace{-.3in}\underbracedmatrixrr{X\\W}{(n-s-l)\,\,\rm columns}\,\,\,\,\right\},
  \tikz[remember picture,overlay]   \draw[dashed,dash pattern={on 4pt off 2pt}] ([xshift=0.5\tabcolsep,yshift=7pt]a.north) -- ([xshift=0.5\tabcolsep,yshift=-2pt]b.south);\tikz[remember picture,overlay]   \draw[dashed,dash pattern={on 4pt off 2pt}] ([xshift=0.5\tabcolsep,yshift=7pt]c.north) -- ([xshift=0.5\tabcolsep,yshift=-2pt]d.south);\tikz[remember picture,overlay]   \draw[dashed,dash pattern={on 4pt off 2pt}] ([xshift=0.5\tabcolsep,yshift=7pt]e.north) -- ([xshift=0.5\tabcolsep,yshift=-2pt]f.south);\tikz[remember picture,overlay]   \draw[dashed,dash pattern={on 4pt off 2pt}] ([xshift=0.5\tabcolsep,yshift=7pt]g.north) -- ([xshift=0.5\tabcolsep,yshift=-2pt]h.south);
\end{equation}
and equipped  with the holomorphic coordinates
\begin{equation}\label{ulx}
    Z:=\left(\begin{matrix}
    z_{11}&\cdots&z_{1(s-p+l)}\\
    \vdots&\ddots&\vdots\\
    z_{l1}&\cdots&z_{l(s-p+l)}\\
    \end{matrix}\right)\,\,,\,\,\,\,X:=\left(\begin{matrix}x_{1(s+l+1)}&\cdots&x_{1n}\\ \vdots&\ddots&\vdots\\ x_{l(s+l+1)}&\cdots &x_{ln}\\
\end{matrix}\right)\,,\,\,\,\,\,\,\,\,\,\,\,\,\,\,\,\,\,\,\,\,\,\,\,\,\,\,\,\,\,
\end{equation}
\begin{equation}
\,\,Y:=\left(\begin{matrix}
y_{(l+1)1}&\cdots& y_{(l+1)(s-p+l)}\\
\vdots&\ddots&\vdots\\
y_{p1}&\cdots& y_{p(s-p+l)}\\
\end{matrix}  \right)\,\,,\,\,\,\,W:=\left(\begin{matrix}
w_{(l+1)(s+l+1)}&\cdots& w_{(l+1)n}\\
\vdots&\ddots&\vdots\\
w_{p(s+l+1)}&\cdots& w_{pn}\\
\end{matrix}  \right) \,.
\end{equation}

Following Sections 3.2 and 3.3 in \cite{F}, we define the following holomorphic atlas $\left\{\left(A^{\tau},(J_l^{\tau})^{-1}\right)\right\}$ for $R_{s,p,n}^{-1}(U_l)$ (the Van der Waerden representation). 

\subsection{The case \texorpdfstring{$p\leq n-s$}{ee} } \label{vanderl} For $0\leq l\leq p$, define an index set $\mathbb J_l$ by \begin{equation}\label{tt}
\mathbb J_l:=\left\{\left(
\begin{matrix}
i_1&i_2&\cdots&i_{p-l}&\cdots&i_p\\
j_1&j_2&\cdots&j_{p-l}&\cdots&j_p\\
\end{matrix}\right)\rule[-.38in]{0.01in}{.82in}\,\,\footnotesize\begin{matrix}
l+1\leq \, i_k\,\leq p\,\,{\rm for}\,\,1\leq\,k\,\leq p-l\,;\,\,\\
1\leq \, i_k\,\leq l\,\,{\rm for}\,\,p-l+1\leq\,k\,\leq p\,;\,\,\\
s+l+1\leq \, j_k\,\leq n\,\,{\rm for}\,\,1\leq\,k\,\leq p-l\,;\,\,\\
1\leq \, j_k\,\leq s-p+l\,\,{\rm for}\,\,p-l+1\leq\,k\,\leq p\,;\,\,\\
i_{k_1}\neq i_{k_2}\,\,{\rm and\,\,}j_{k_1}\neq j_{k_2}\,\,{\rm for\,\,} k_1\neq k_2.\end{matrix}\right\}.
\end{equation}
Associate each $\tau=\left(\begin{matrix}
i_1&i_2&\cdots&i_p\\
j_1&j_2&\cdots&j_p\\
\end{matrix}\right)\in\mathbb J_l$ with a complex Euclidean space $\mathbb {C}^{p(n-p)}$ equipped with holomorphic coordinates  $\left(\widetilde X,\widetilde Y,\overrightarrow B^1,\cdots,\overrightarrow B^p\right)$ defined as follows.
\begin{equation}\label{ulu}
\widetilde X:=\left(\begin{matrix}
x_{1(s+l+1)}&\cdots &x_{1n}\\
\vdots&\ddots&\vdots\\
x_{l(s+l+1)}&\cdots &x_{ln}\\
\end{matrix}\right)\,\,\,{\rm and}\,\,\,\widetilde Y:=\left(\begin{matrix}
y_{(l+1)1}&\cdots& y_{(l+1)(s-p+l)}\\
\vdots&\ddots&\vdots\\
y_{p1}&\cdots& y_{p(s-p+l)}\\
\end{matrix}  \right)  \,;\,
\end{equation}
for $1\leq k\leq p-l$,
\begin{equation}\label{qbp}
\begin{split}
&\overrightarrow B^{k}:=\left(b_{i_{k}j_{k}},\xi^{(k)}_{i_{k}(s+l+1)},\xi^{(k)}_{i_{k}(s+l+2)},\cdots,\widehat{\xi^{(k)}_{i_{k}j_1}},\cdots,\widehat{\xi^{(k)}_{i_{k}j_2}},\cdots,\widehat{\xi^{(k)}_{i_{k}j_{k}}},\cdots,\xi^{(k)}_{i_{k}n},\right.\\
&\,\,\,\,\,\,\,\,\,\,\,\,\,\,\,\,\,\,\,\,\,\,\,\,\,\,\,\,\,\,\left.\xi^{(k)}_{(l+1)j_{k}},\xi^{(k)}_{(l+2)j_{k}},\cdots,\widehat{\xi^{(k)}_{i_1j_{k}}},\cdots,\widehat{\xi^{(k)}_{i_2j_{k}}},\cdots,\widehat{\xi^{(k)}_{i_{k}j_{k}}},\cdots,\xi^{(k)}_{pj_{k}}\right)\,;\\
\end{split}
\end{equation}
for $p-l+1\leq k\leq p$,
\begin{equation}\label{qbpb}
\begin{split}
&\overrightarrow B^{k}:=\left(a_{i_{k}j_{k}},\xi^{(k)}_{i_{k}1},\xi^{(k)}_{i_{k}2},\cdots,\widehat{\xi^{(k)}_{i_{k}j_{p-l+1}}},\cdots,\widehat{\xi^{(k)}_{i_{k}j_{p-l+2}}},\cdots,\widehat{\xi^{(k)}_{i_{k}j_{k}}},\cdots,\xi^{(k)}_{i_{k}(s-p+l)},\right.\\
&\,\,\,\,\,\,\,\,\,\,\,\,\,\,\,\,\,\,\,\,\,\,\,\,\,\,\,\,\,\,\left.\xi^{(k)}_{1j_{k}},\xi^{(k)}_{2j_{k}},\cdots,\widehat{\xi^{(k)}_{i_{p-l+1}j_{k}}},\cdots,\widehat{\xi^{(k)}_{i_{p-l+2}j_{k}}},\cdots,\cdots,\widehat{\xi^{(k)}_{i_{k}j_{k}}},\cdots,\xi^{(k)}_{lj_{k}}\right)\,.\,\,\,\,\,\,\,\,\,\,\,\,\,\\
\end{split}
\end{equation}

The holomorphic embedding  $J_l^{\tau}:\mathbb C^{p(n-p)}\rightarrow\mathcal T_{s,p,n}\hookrightarrow\mathbb {CP}^{N_{p,n}}\times\mathbb {CP}^{N^0_{s,p,n}}\times\cdots\times\mathbb {CP}^{N^p_{s,p,n}}$ is the holomorphic extension of the birational map $\mathcal K_{s,p,n}\circ \Gamma_l^{\tau}$, 
where $\mathcal K_{s,p,n}$ is given by (\ref{jl}) and $\Gamma_l^{\tau}:\mathbb C^{p(n-p)}\rightarrow U_l$ is defined by
\begin{equation}\label{ws}
\begin{split}
&\,\,\,\,\,\,\,\,\,\,\,\,\,\,\,\,\,\,\,\,\,\,\,\,\,\,\,\,\,\,\,\,\,\,\,\,\,\,\,\,\,\,\,\,\,\,\,\,\,\,\,\,\,\,\,\,\,\,\,\,\,\,\,\,\,\Gamma_l^{\tau}\left(\widetilde X,\widetilde Y,\overrightarrow B^1,\cdots,\overrightarrow B^p\right):=\\
&\left(
\begin{matrix}
\sum\limits_{k=p-l+1}^p\left(\prod\limits_{t=p-l+1}^{k}a_{i_{t}j_t}\right)\cdot\Xi_k^T\cdot\Omega_k &0_{l\times(p-l)}&I_{l\times l}&\widetilde X\\ \widetilde Y&I_{(p-l)\times(p-l)}&0_{(p-l)\times l}&\sum\limits_{k=1}^{p-l}\left(\prod\limits_{t=1}^{k}b_{i_tj_t}\right)\cdot\Xi_k^T\cdot\Omega_k\\
\end{matrix}\right)\,.      \end{split}
\end{equation}
Here $\Xi_k$ and $\Omega_k$ are defined as follows. For $1\leq k\leq p-l$, $\Xi_k:=
\left(v_{l+1}^k,\cdots,v_p^k\right)$ where
\begin{equation}\label{w3}
     v_t^k=\left\{\begin{array}{ll}
    \xi^{(k)}_{tj_k} \,\,\,\,\,\,\,\,\,\,\,\,\,\,\,\,\,\,\,\,\,\,\,\,\,\,\,\,\,\, & t\in\{l+1,l+2,\cdots,p\}\backslash\{i_1,i_2,\cdots,i_k\}\,\,\,\,\,\,\,\,\,\,\,\,\,\,\,\,\, \,\,\,\,\,\,\,\,\,\,\,\,\,\,\,\,\,\\
    0& t\in\{i_1,i_2,\cdots,i_{k-1}\}\,\,\,\,\,\,\, \\
  
    1 \,\,\,\,\,\,\,\,\,\,\,\,\,\,\,\,\,\,\,\,\,\,\,\,\,\,\,\,\,\, &t=i_k\,\,\,\,\,\,\,\,\,\,\,\,\,\,\,\,\, \,\,\,\,\,\,\,\,\,\,\,\,\,\,\,\,\,\\
    \end{array}\right.,
\end{equation}
and $\Omega_k:=
\left(w_{s+l+1}^k,\cdots,w_n^k\right)$ where
\begin{equation}\label{w4}
     w_t^k=\left\{\begin{array}{ll}
    \xi^{(k)}_{i_kt} \,\,\,\,\,\,\,\,\,\,\,\,\,\,\,\,\,\,\,\,\,\,\,\,\,\,\,\,\,\, & t\in\{s+l+1,s+l+2,\cdots,n\}\backslash\{j_1,j_2,\cdots,j_k\}\,\,\\
    0& t\in\{j_1,j_2,\cdots,j_{k-1}\}\,\,\,\,\,\,\, \\
  
    1 \,\,\,\,\,\,\,\,\,\,\,\,\,\,\,\,\,\,\,\,\,\,\,\,\,\,\,\,\,\, &t=j_k\,\,\,\,\,\,\,\,\,\\
    \end{array}\right..
\end{equation}
For $p-l+1\leq k\leq p$,  $\Xi_k:=
\left(v_{1}^k,\cdots,v_l^k\right)$ where
\begin{equation}\label{w5}
     v_t^k=\left\{\begin{array}{ll}
    \xi^{(k)}_{tj_k} \,\,\,\,\,\,\,\,\,\,\,\,\,\,\,\,\,\,\,\,\,\,\,\,\,\,\,\,\,\, & t\in\{1,2,\cdots,l\}\backslash\{i_{p-l+1},i_{p-l+2},\cdots,i_k\}\,\,\,\,\,\,\,\,\,\,\,\,\,\,\,\,\, \,\,\,\,\,\,\,\,\,\,\,\,\,\,\,\,\,\\
    0& t\in\{i_{p-l+1},i_{p-l+2},\cdots,i_{k-1}\}\,\,\,\,\,\,\, \\
  
    1 \,\,\,\,\,\,\,\,\,\,\,\,\,\,\,\,\,\,\,\,\,\,\,\,\,\,\,\,\,\, &t=i_k\,\,\,\,\,\,\,\,\,\,\,\,\,\,\,\,\, \,\,\,\,\,\,\,\,\,\,\,\,\,\,\,\,\,\\
    \end{array}\right.,
\end{equation}
and $\Omega_k:=
\left(w_{1}^k,\cdots,w_{s-p+l}^k\right)$ where
\begin{equation}\label{w2}
     w_t^k=\left\{\begin{array}{ll}
    \xi^{(k)}_{i_kt} \,\,\,\,\,\,\,\,\,\,\,\,\,\,\,\,\,\,\,\,\,\,\,\,\,\,\,\,\,\, & t\in\{1,2,\cdots,s-p+l\}\backslash\{j_{p-l+1},j_{p-l+2},\cdots,j_k\} \,\,\,\,\,\,\,\,\,\,\,\,\,\,\,\,\,\\
    0& t\in\{j_{p-l+1},j_{p-l+2},\cdots,j_{k-1}\}\,\,\,\,\,\,\, \\
  
    1 \,\,\,\,\,\,\,\,\,\,\,\,\,\,\,\,\,\,\,\,\,\,\,\,\,\,\,\,\,\, &t=j_k\,\,\,\,\,\,\,\,\,\,\,\,\,\\
    \end{array}\right..
\end{equation}

Let $A^{\tau}$ be the image of $\mathbb C^{p(n-p)}$ under the holomorphic map $J_l^{\tau}$.  Then, $\left\{\left(A^{\tau},(J_l^{\tau})^{-1}\right)\right\}_{\tau\in\mathbb J_p}$  is a holomorphic atlas of $R^{-1}_{s,p,n}(U_l)$.

\begin{example}
Consider Grassmannian $G(4,8)$, with $s=2$, $l=2$, and $\tau=\left(
\begin{matrix}
3&4&1&2\\
7&8&1&2\\
\end{matrix}\right)\in\mathbb J_2$. The holomorphic coordinates $\left(\widetilde X,\widetilde Y,\overrightarrow B^1,\overrightarrow B^2, \overrightarrow B^3, \overrightarrow B^4\right)$ of $\mathbb {C}^{16}$ are	\begin{equation}
\begin{split}
&\widetilde X=\left(x_{17},x_{18},x_{27},x_{28}\right),\,\,\,\,\,\,\widetilde Y=\left(y_{31},y_{32},y_{41},y_{42}\right),\\
&\overrightarrow  B^1=\left(b_{37},\xi^{(1)}_{38},\xi^{(1)}_{47}\right),\,\,\,\,\,\,\overrightarrow  B^2=\left(b_{48}\right),\,\,\,\,\,\, B^3=\left(a_{11},\xi^{(3)}_{12},\xi^{(3)}_{21}\right),\,\,\,\,\,\,\overrightarrow  B^4=\left(a_{22}\right).\\
\end{split}
\end{equation}
The holomorphic map $\Gamma_2^{\tau}:\mathbb C^{16}\rightarrow U_2$ is given by $\Gamma^{\tau}\left(\widetilde X,\widetilde Y,\overrightarrow B^1,\cdots,\overrightarrow B^4\right)=$
\begin{equation}
\left(
\begin{matrix}
a_{11}&a_{1}\cdot\xi^{(3)}_{12}&0&0&1&0&x_{17}&x_{18}\\
a_{11}\cdot\xi^{(3)}_{21}&a_{11}\cdot(\xi^{(3)}_{12}\cdot\xi^{(3)}_{21}+a_{22})&0&0&0&1&x_{27}&x_{28}\\
y_{31}&y_{32}&1&0&0&0&b_{37}&b_{37}\cdot\xi^{(1)}_{38}\\
y_{41}&y_{42}&0&1&0&0&b_{37}\cdot\xi^{(1)}_{47}&b_{37}\cdot(\xi^{(1)}_{47}\cdot\xi^{(1)}_{38}+b_{48})\\
	\end{matrix}
	\right).
	\end{equation}
\end{example}

\subsection{The case \texorpdfstring{$n-s<p$}{ff}} \label{vandern} For $0\leq l\leq n-s$, define an index set $\mathbb J_l$ by
	\begin{equation}
 \left\{\left(
	\begin{matrix}
	i_1&\cdots&i_{n-s}\\
	j_1&\cdots&j_{n-s}\\
	\end{matrix}\right)\rule[-.35in]{0.01in}{.72in}\footnotesize\begin{matrix}
	(i_{n-s-l+1},\,i_{n-s-l+2},\,\cdots,\,i_{n-s})\,\, {\rm is\,\,a\,\,permutation\,\,of\,\,}(1,\,2,\,\cdots,\,\,l)\,;\\
	(j_1,\,\cdots,\,j_{n-s-l})\,\,{\rm is\,\,a\,\,permutations\,\,of\,\,}(s+l+1,\,s+l+2,\,\cdots,\,n)\,;\\
	1\leq\, j_t\,\leq s-p+l\,\,{\rm for}\,\,n-s-l+1\leq\,t\,\leq n-s\,,\,\,j_{t_1}\neq j_{t_2}\,\,{\rm for\,\,} t_1\neq t_2\,;\\
	l+1\leq\, i_t\,\leq p\,\,{\rm for}\,\,1\leq\,t\,\leq n-s-l\,,i_{t_1}\neq i_{t_2}\,\,{\rm for\,\,} t_1\neq t_2\\
	\end{matrix}\right\}.
	\end{equation}
Associate each $\tau=\left(\begin{matrix}
i_1&i_2&\cdots&i_{n-s}\\
j_1&j_2&\cdots&j_{n-s}\\
\end{matrix}\right)\in\mathbb J_l$ with a complex Euclidean space $\mathbb {C}^{p(n-p)}$ equipped with the holomorphic coordinates  $\left(\widetilde X,\widetilde Y,\overrightarrow B^1,\cdots,\overrightarrow B^{n-s}\right)$ defined as follows.
\begin{equation}\label{rulu}
\widetilde X:=\left(\begin{matrix}
x_{1(s+l+1)}&\cdots &x_{1n}\\
\vdots&\ddots&\vdots\\
x_{l(s+l+1)}&\cdots &x_{ln}\\
\end{matrix}\right)\,\,\,\,{\rm and}\,\,\,\,\,\widetilde Y:=\left(\begin{matrix}
y_{(l+1)1}&\cdots& y_{(l+1)(s-p+l)}\\
\vdots&\ddots&\vdots\\
y_{p1}&\cdots& y_{p(s-p+l)}\\
\end{matrix}  \right)\,;\,\,
\end{equation}
for $1\leq k\leq n-s-l$,
\begin{equation}\label{rqbp}
\begin{split}
&\overrightarrow B^{k}:=\left(a_{i_{k}j_{k}},\xi^{(k)}_{i_{k}(s+l+1)},\xi^{(k)}_{i_{k}(s+l+2)},\cdots,\widehat{\xi^{(k)}_{i_{k}j_1}},\cdots,\widehat{\xi^{(k)}_{i_{k}j_2}},\cdots,\widehat{\xi^{(k)}_{i_{k}j_{k}}},\cdots,\xi^{(k)}_{i_{k}n},\right.\\
&\,\,\,\,\,\,\,\,\,\,\,\,\,\,\,\,\,\,\,\,\,\,\,\,\,\,\,\,\left.\xi^{(k)}_{(l+1)j_{k}},\xi^{(k)}_{(l+2)j_{k}},\cdots,\widehat{\xi^{(k)}_{i_1j_{k}}},\cdots,\widehat{\xi^{(k)}_{i_2j_{k}}},\cdots,\cdots,\widehat{\xi^{(k)}_{i_{k}j_{k}}},\cdots,\xi^{(k)}_{pj_{k}}\right)\,;\\
\end{split}
\end{equation}
for $n-s-l+1\leq k\leq n-s$,
\begin{equation}\label{rqbpb}
\begin{split}
&\overrightarrow B^{k}:=\left(a_{i_{k}j_{k}},\xi^{(k)}_{i_{k}1},\xi^{(k)}_{i_{k}2},\cdots,\widehat{\xi^{(k)}_{i_{k}j_{n-s-l+1}}},\cdots,\widehat{\xi^{(k)}_{i_{k}j_{n-s-l+2}}},\cdots,\widehat{\xi^{(k)}_{i_{k}j_{k}}},\cdots,\xi^{(k)}_{i_{k}(s-p+l)},\right.\\
&\,\,\,\,\,\,\,\,\,\,\,\,\,\,\,\,\,\,\,\,\,\,\,\,\,\,\,\,\left.\xi^{(k)}_{1j_{k}},\xi^{(k)}_{2j_{k}},\cdots,\widehat{\xi^{(k)}_{i_{n-s-l+1}j_{k}}},\cdots,\widehat{\xi^{(k)}_{i_{n-s-l+2}j_{k}}},\cdots,\cdots,\widehat{\xi^{(k)}_{i_{k}j_{k}}},\cdots,\xi^{(k)}_{lj_{k}}\right)\,.\,\,\,\,\,\,\,\,\,\,\,\,\,\,\\
\end{split}
\end{equation}

The holomorphic embedding $J_l^{\tau}:\mathbb C^{p(n-p)}\rightarrow\mathcal T_{s,p,n}\hookrightarrow\mathbb {CP}^{N_{p,n}}\times\mathbb {CP}^{N^0_{s,p,n}}\times\cdots\times\mathbb {CP}^{N^p_{s,p,n}}$ is the holomorphic extension of the birational map $\mathcal K_{s,p,n}\circ \Gamma_l^{\tau}$, where $\Gamma_l^{\tau}:\mathbb C^{p(n-p)}\rightarrow U_l$ is given by
\begin{equation}\label{gammalr}
\small
\begin{split}
&\,\,\,\,\,\,\,\,\,\,\,\,\,\,\,\,\,\,\,\,\,\,\,\,\,\,\,\,\,\,\,\,\,\,\,\,\,\,\,\,\,\,\,\,\,\,\,\,\,\,\,\,\,\,\,\,\,\,\,\,\,\,\,\,\,\Gamma_l^{\tau}\left(\widetilde X,\widetilde Y,\overrightarrow B^1,\cdots,\overrightarrow B^{n-s}\right):=\\
&\left(
\begin{matrix}
\sum\limits_{k=n-s-l+1}^{n-s}\left(\prod\limits_{t=n-s-l+1}^{k}a_{i_{t}j_t}\right)\cdot\Xi_k^T\cdot\Omega_k &0_{l\times(p-l)}&I_{l\times l}&\widetilde X\\ \widetilde Y&I_{(p-l)\times(p-l)}&0_{(p-l)\times l}&\sum\limits_{k=1}^{n-s-l}\left(\prod\limits_{t=1}^{k}b_{i_tj_t}\right)\cdot\Xi_k^T\cdot\Omega_k\\
\end{matrix}\right)\,.       
\end{split}
\end{equation}
For $1\leq k\leq n-s-l$, $\Xi_k:=
\left(v_{l+1}^k,\cdots,v_p^k\right)$ where
\begin{equation}\label{lw1}
     v_t^k=\left\{\begin{array}{ll}
    \xi^{(k)}_{tj_k} \,\,\,\,\,\,\,\,\,\,\,\,\,\,\,\,\,\,\,\,\,\,\,\,\,\,\,\,\,\, & t\in\{l+1,l+2,\cdots,p\}\backslash\{i_1,i_2,\cdots,i_k\}\,\,\,\,\,\,\,\,\,\,\,\,\,\,\,\,\, \,\,\,\,\,\,\,\,\,\,\,\,\,\,\,\,\,\\
    0& t\in\{i_1,i_2,\cdots,i_{k-1}\}\,\,\,\,\,\,\, \\
  
    1 \,\,\,\,\,\,\,\,\,\,\,\,\,\,\,\,\,\,\,\,\,\,\,\,\,\,\,\,\,\, &t=i_k\,\,\,\,\,\,\,\,\,\,\,\,\,\,\,\,\, \,\,\,\,\,\,\,\,\,\,\,\,\,\,\,\,\,\\
    \end{array}\right.,
\end{equation}
and $\Omega_k:=
\left(w_{s+l+1}^k,\cdots,w_n^k\right)$ where
\begin{equation}\label{lw2}
     w_t^k=\left\{\begin{array}{ll}
    \xi^{(k)}_{i_kt} \,\,\,\,\,\,\,\,\,\,\,\,\,\,\,\,\,\,\,\,\,\,\,\,\,\,\,\,\,\, & t\in\{s+l+1,s+l+2,\cdots,n\}\backslash\{j_1,j_2,\cdots,j_k\}\,\,\\
    0& t\in\{j_1,j_2,\cdots,j_{k-1}\}\,\,\,\,\,\,\, \\
  
    1 \,\,\,\,\,\,\,\,\,\,\,\,\,\,\,\,\,\,\,\,\,\,\,\,\,\,\,\,\,\, &t=j_k\,\,\,\,\,\,\,\,\,\\
    \end{array}\right..
\end{equation}
For $n-s-l+1\leq k\leq n-s$, $\Xi_k:=
\left(v_{1}^k,\cdots,v_l^k\right)$ where
\begin{equation}\label{lw3}
     v_t^k=\left\{\begin{array}{ll}
    \xi^{(k)}_{tj_k} \,\,\,\,\,\,\,\,\,\,\,\,\,\,\,\,\,\,\,\,\,\,\,\,\,\,\,\,\,\, & t\in\{1,2,\cdots,l\}\backslash\{i_{n-s-l+1},i_{n-s-l+2},\cdots,i_k\}\,\,\,\,\,\,\,\,\,\,\,\,\,\,\,\,\, \,\,\,\,\,\,\,\,\,\,\,\,\,\,\,\,\,\\
    0& t\in\{i_{n-s-l+1},i_{n-s-l+2},\cdots,i_{k-1}\}\,\,\,\,\,\,\, \\
  
    1 \,\,\,\,\,\,\,\,\,\,\,\,\,\,\,\,\,\,\,\,\,\,\,\,\,\,\,\,\,\, &t=i_k\,\,\,\,\,\,\,\,\,\,\,\,\,\,\,\,\, \,\,\,\,\,\,\,\,\,\,\,\,\,\,\,\,\,\\
    \end{array}\right.,
\end{equation}
and $\Omega_k:=
\left(w_{1}^k,\cdots,w_{s-p+l}^k\right)$ where
\begin{equation}\label{lw4}
     w_t^k=\left\{\begin{array}{ll}
    \xi^{(k)}_{i_kt} \,\,\,\,\,\,\,\,\,\,\,\,\,\,\,\,\,\,\,\,\,\,\,\,\,\,\,\,\,\, & t\in\{1,2,\cdots,s-p+l\}\backslash\{j_{n-s-l+1},j_{n-s-l+2},\cdots,j_k\} \,\,\,\,\,\,\,\,\,\,\,\,\,\,\,\,\,\\
    0& t\in\{j_{n-s-l+1},j_{n-s-l+2},\cdots,j_{k-1}\}\,\,\,\,\,\,\, \\
  
    1 \,\,\,\,\,\,\,\,\,\,\,\,\,\,\,\,\,\,\,\,\,\,\,\,\,\,\,\,\,\, &t=j_k\,\,\,\,\,\,\,\,\,\,\,\,\,\\
    \end{array}\right..
\end{equation}

Let $A^{\tau}$ be the image of $\mathbb C^{p(n-p)}$ under $J_l^{\tau}$. Then, $\left\{\left(A^{\tau},(J^{\tau})^{-1}\right)\right\}_{\tau\in\mathbb J_l}$  is a holomorphic atlas for $R^{-1}_{s,p,n}(U_l)$, $0\leq l\leq n-s$.

\section{Holomorphic Atlas}\label{section:ncc}

\subsection{Proof of Lemma \ref{rwk2}.} \label{dvanderl} 
Since $D^+_{1},D^+_2,\cdots,D^+_{r+1-j}$ have empty intersection with $D^-_j$, by Lemma \ref{bb} we can conclude the following formulas for line bundles.
\begin{equation}\label{rdbst-1}
\begin{split}
\check B^{-j}_0&=(R_{s,p,n})^*\left(\mathcal O_{G(p,n)}(1)\right)|_{D^-_j}-\sum_{i=r+2-j}^{r}(r+1-i)\cdot D^{-j}_{+i}\,,\\
\check B^{-j}_1&=(R_{s,p,n})^*\left(\mathcal O_{G(p,n)}(1)\right)|_{D^-_j}-\sum_{i=r+2-j}^{r-1}(r-i)\cdot D^{-j}_{+i}-D^{-j}_{-1}\,,\\
\check B^{-j}_2&=(R_{s,p,n})^*\left(\mathcal O_{G(p,n)}(1)\right)|_{D^-_j}-\sum_{i=r+2-j}^{r-2}(r-1-i)\cdot D^{-j}_{+i}-2D^{-j}_{-1}-D^{-j}_{-2}\,,\\
&\,\,\,\vdots\\
\end{split}
\end{equation}
\begin{equation}\label{rdm-2}
\begin{split}
\check B^{-j}_{j-2}&=(R_{s,p,n})^*\left(\mathcal O_{G(p,n)}(1)\right)|_{D^-_j}-D^{-j}_{+(r+2-j)}-(j-2)D^{-j}_{-1}-(j-3)D^{-j}_{-2}-\cdots-D^{-j}_{-(j-2)}\,,\\
\end{split}
\end{equation}
\begin{equation}\label{rdm-1}
\begin{split}
\check B^{-j}_{j-1}&=(R_{s,p,n})^*\left(\mathcal O_{G(p,n)}(1)\right)|_{D^-_j}-(j-1)D^{-j}_{-1}-(j-2)D^{-j}_{-2}-\cdots-D^{-j}_{-(j-1)}\,,\\
\end{split}
\end{equation}
\begin{equation}\label{rdm}
\begin{split}
\check B^{-j}_j&=(R_{s,p,n})^*\left(\mathcal O_{G(p,n)}(1)\right)|_{D^-_j}-jD^{-j}_{-1}-(j-1)D^{-j}_{-2}-\cdots-2D^{-j}_{-(j-1)}-D^{-j}_{-j}\,,\\
\end{split}
\end{equation}
\begin{equation}\label{rdm+1}
\begin{split}
\check B^{-j}_{j+1}&=(R_{s,p,n})^*\left(\mathcal O_{G(p,n)}(1)\right)|_{D^-_j}-(j+1)D^{-j}_{-1}-jD^{-j}_{-2}-\cdots-2D^{-j}_{-j}-D^{-j}_{-(j+1)}\,,\\
\end{split}
\end{equation}
\begin{equation}
\begin{split}
&\,\,\,\vdots\\
\check B^{-j}_r&=D^{-j}_{-r}=(R_{s,p,n})^*\left(\mathcal O_{G(p,n)}(1)\right)|_{D^-_j}-rD^{-j}_{-1}-\cdots-\cdots-2D^{-j}_{-(r-1)}\,.\\
\end{split}
\end{equation} 
Then (\ref{1rd-j-1}) and (\ref{1rd-j1-1}) follow from (\ref{rdm-2}), (\ref{rdm-1}), (\ref{rdm}), and (\ref{rdm+1}) directly.

When $1\leq j\leq r-2$, subtracting (\ref{rdm+1}) from (\ref{rdm}), we derive that
\begin{equation}\label{rdjj}
D^{-j}_{-1}+D^{-j}_{-2}+\cdots+D^{-j}_{-(j-1)}+D^{-j}_{-j}+D^{-j}_{-(j+1)}=\check B^{-j}_j-\check B^{-j}_{j+1}.\\
\end{equation}
Restricting (\ref{tk2}) to $D^-_j$ and plugging in (\ref{rdjj}), we have that
\begin{equation}
\begin{split}
 -K_{\mathcal T_{s,p,n}}|_{D^{- }_j}=&(s-p+1)\check B^{-j}_0 +2\sum_{m=1}^{r-1}\check B^{-j}_m +\sum_{i=1}^{j+1}D^{-j}_{-i}+
\sum_{i=j+2}^rD^{-j}_{-i}+\sum_{i=r+2-j}^rD^{-j}_{+i}\,\\
=&(s-p+1) \check B^{-j}_0 +2\sum_{m=1}^{r-1}\check B^{-j}_m +\check B^{-j}_j-\check B^{-j}_{j+1}+\sum_{i=j+2}^{r}D^{-j}_{-i}+\sum_{i=r+2-j}^rD^{-j}_{+i}\,.\\
\end{split}
\end{equation}
Similarly, when $j=r-1$ or $r$, we have that
\begin{equation}
D^{-j}_{-1}+D^{-j}_{-2}+\cdots+D^{-j}_{-(j-1)}+D^{-j}_{-j}+D^{-j}_{-(j+1)}=\check B^{-j}_{r-1},\\
\end{equation}
where we make the convention that $D^{-r}_{-(r+1)}$ is trivial; hence,
\begin{equation}
\begin{split}
 -K_{\mathcal T_{s,p,n}}|_{D^{- }_j}=&(s-p+1)\check B^{-j}_0 +2\sum_{m=1}^{r-1}\check B^{-j}_m +\sum_{i=1}^{j+1}D^{-j}_{-i}+\sum_{i=j+2}^rD^{-j}_{-i}+\sum_{i=r+2-j}^rD^{-j}_{+i}\,\\
=&(s-p+1) \check B^{-j}_0 +2\sum_{m=1}^{r-1}\check B^{-j}_m +\check B^{-j}_{r-1}+\sum_{i=j+2}^{r}D^{-j}_{-i}+\sum_{i=r+2-j}^rD^{-j}_{+i}\,.\\
\end{split}
\end{equation}

Next we consider the restriction on $D^{+}_j$. By Lemma \ref{bb} and the fact that $D^-_{1},D^-_2,\cdots,D^-_{r+1-j}$ have empty intersection with $D^+_j$, we can conclude that
\begin{equation}\label{r+dbst}
\begin{split}
\check B^{+j}_0&=(R_{s,p,n})^*\left(\mathcal O_{G(p,n)}(1)\right)|_{D^+_j}-\sum_{i=1}^{r}(r+1-i)\cdot D^{+j}_{+i}\,,\\
\check B^{+j}_1&=(R_{s,p,n})^*\left(\mathcal O_{G(p,n)}(1)\right)|_{D^+_j}-\sum_{i=1}^{r-1}(r-i)\cdot D^{+j}_{+i}\,,\\
&\,\,\,\vdots\\
\end{split}
\end{equation}
\begin{equation}\label{r+dm-2}
\begin{split}
\check B^{+j}_{r-1-j}&=(R_{s,p,n})^*\left(\mathcal O_{G(p,n)}(1)\right)|_{D^+_j}-\sum_{i=1}^{j+1}(j+2-i)\cdot D^{+j}_{+i}\,,\\
\end{split}
\end{equation}
\begin{equation}\label{r+dm-1}
\begin{split}
\check B^{+j}_{r-j}&=(R_{s,p,n})^*\left(\mathcal O_{G(p,n)}(1)\right)|_{D^+_j}-\sum_{i=1}^{j}(j+1-i)\cdot D^{+j}_{+i}\,,\\
\end{split}
\end{equation}
\begin{equation}\label{r+dm}
\begin{split}
\check B^{+j}_{r+1-j}&=(R_{s,p,n})^*\left(\mathcal O_{G(p,n)}(1)\right)|_{D^+_j}-\sum_{i=1}^{j-1}(j-i)\cdot D^{+j}_{+i}\,,\\
\end{split}
\end{equation}
\begin{equation}\label{r+dm+1}
\begin{split}
\check B^{+j}_{r+2-j}&=(R_{s,p,n})^*\left(\mathcal O_{G(p,n)}(1)\right)|_{D^+_j}-\sum_{i=1}^{j-2}(j-1-i)\cdot D^{+j}_{+i}-D^{+j}_{-(r+2-j)}\,,\\
\end{split}
\end{equation}
\begin{equation}\label{116}
\begin{split}
&\,\,\,\vdots\\
\check B^{+j}_{r-1}&=(R_{s,p,n})^*\left(\mathcal O_{G(p,n)}(1)\right)|_{D^+_j}-D^{+j}_{+1}-\sum_{i=r+2-j}^{r-1}(r-i)\cdot D^{+j}_{-i}\,.\\
\check B^{+j}_r=D^{+j}_{-r}&=(R_{s,p,n})^*\left(\mathcal O_{G(p,n)}(1)\right)|_{D^+_j}-\sum_{i=r+2-j}^{r-1}(r+1-i)\cdot D^{+j}_{-i}\,.\\
\end{split}
\end{equation} 

Then (\ref{r+d-j-1}) and (\ref{r+d-j1-1}) follow from (\ref{r+dm-2}), (\ref{r+dm-1}), (\ref{r+dm}),  (\ref{r+dm+1}), and (\ref{116}).

When $1\leq j\leq r-1$, subtracting (\ref{r+dm-2}) from (\ref{r+dm-1}), we derive that
\begin{equation}\label{r+djj}
D^{+j}_{+1}+D^{+j}_{+2}+\cdots+D^{+j}_{+(j-1)}+D^{+j}_{+j}+D^{+j}_{+(j+1)}=\check B^{+j}_{r-j}-\check B^{+j}_{r-1-j}.\\
\end{equation}
Restricting (\ref{tk2}) to $D^+_j$ and plugging in (\ref{r+djj}), we have that
\begin{equation}
\begin{split}
 -K_{\mathcal T_{s,p,n}}|_{D^{+}_j}&=(s-p+1)\cdot \check B^{+j}_0 +2\sum_{m=1}^{r-1}\check B^{+j}_m +\sum_{i=1}^{j+1}D^{+j}_{+i}+
\sum_{i=j+2}^rD^{+j}_{+i}+\sum_{i=r+2-j}^rD^{+j}_{-i}\,,\\
=&(s-p+1)\cdot \check B^{+j}_0 +2\sum_{m=1}^{r-1}\check B^{+j}_m +\check B^{+j}_{r-j}-\check B^{+j}_{r-1-j}+\sum_{i=j+2}^{r}D^{+j}_{+i}+\sum_{i=r+2-j}^rD^{+j}_{-i}\,.\\
\end{split}
\end{equation}

When $j=r$, subtracting (\ref{r+dm-2}) from (\ref{r+dm-1}), we derive that
\begin{equation}\label{r+djj1}
D^{+j}_{+1}+D^{+j}_{+2}+\cdots+D^{+j}_{+(j-1)}+D^{+j}_{+j}+D^{+j}_{+(j+1)}=\check B^{+j}_{1}-\check B^{+j}_{0}.\\
\end{equation}
Restricting (\ref{tk2}) to $D^+_j$ and plugging in (\ref{r+djj1}), we have that
\begin{equation}
\begin{split}
 -K_{\mathcal T_{s,p,n}}|_{D^{+}_j}&=(s-p+1)\cdot \check B^{+j}_0 +2\sum_{m=1}^{r-1}\check B^{+j}_m +\sum_{i=1}^{j+1}D^{+j}_{+i}+
\sum_{i=j+2}^rD^{+j}_{+i}+\sum_{i=r+2-j}^rD^{+j}_{-i}\,,\\
=&(s-p+1)\cdot \check B^{+j}_0 +2\sum_{m=1}^{r-1}\check B^{+j}_m +\check B^{+j}_{1}-\check B^{+j}_{0}+\sum_{i=j+2}^{r}D^{+j}_{+i}+\sum_{i=r+2-j}^rD^{+j}_{-i}\,.\\
\end{split}
\end{equation}

We complete the proof of Lemma \ref{rwk2}.\,\,\,\,$\endpf$ 

\subsection{Proof of Lemma \ref{rwk3}.}\label{dvanderld}
Since $D^+_{1},D^+_2,\cdots,D^+_{r+1-j}$ have empty intersection with $D^-_j$, by Lemma \ref{bb} we can conclude that
\begin{equation}\label{prdbst-1}
\begin{split}
\check B^{-j}_0&=D^{-j}_{+r}=(R_{s,p,n})^*\left(\mathcal O_{G(p,n)}(1)\right)|_{D^-_j}-\sum_{i=r+2-j}^{r-1}(r+1-i)\cdot D^{-j}_{+i}\,,\\
\check B^{-j}_1&=(R_{s,p,n})^*\left(\mathcal O_{G(p,n)}(1)\right)|_{D^-_j}-\sum_{i=r+2-j}^{r-1}(r-i)\cdot D^{-j}_{+i}-D^{-j}_{-1}\,,\\
\check B^{-j}_2&=(R_{s,p,n})^*\left(\mathcal O_{G(p,n)}(1)\right)|_{D^-_j}-\sum_{i=r+2-j}^{r-2}(r-1-i)\cdot D^{-j}_{+i}-2D^{-j}_{-1}-D^{-j}_{-2}\,,\\
&\,\,\,\vdots\\
\end{split}
\end{equation}
\begin{equation}\label{prdm-2}
\begin{split}
\check B^{-j}_{j-2}&=(R_{s,p,n})^*\left(\mathcal O_{G(p,n)}(1)\right)|_{D^-_j}-D^{-j}_{+(r+2-j)}-(j-2)D^{-j}_{-1}-(j-3)D^{-j}_{-2}-\cdots-D^{-j}_{-(j-2)}\,,\\
\end{split}
\end{equation}
\begin{equation}\label{prdm-1}
\begin{split}
\check B^{-j}_{j-1}&=(R_{s,p,n})^*\left(\mathcal O_{G(p,n)}(1)\right)|_{D^-_j}-(j-1)D^{-j}_{-1}-(j-2)D^{-j}_{-2}-\cdots-D^{-j}_{-(j-1)}\,,\\
\end{split}
\end{equation}
\begin{equation}\label{prdm}
\begin{split}
\check B^{-j}_j&=(R_{s,p,n})^*\left(\mathcal O_{G(p,n)}(1)\right)|_{D^-_j}-jD^{-j}_{-1}-(j-1)D^{-j}_{-2}-\cdots-2D^{-j}_{-(j-1)}-D^{-j}_{-j}\,,\\
\end{split}
\end{equation}
\begin{equation}\label{prdm+1}
\begin{split}
\check B^{-j}_{j+1}&=(R_{s,p,n})^*\left(\mathcal O_{G(p,n)}(1)\right)|_{D^-_j}-(j+1)D^{-j}_{-1}-jD^{-j}_{-2}-\cdots-2D^{-j}_{-j}-D^{-j}_{-(j+1)}\,,\\
\end{split}
\end{equation}
\begin{equation}
\begin{split}
&\,\,\,\vdots\\
\check B^{-j}_{r-1}&=(R_{s,p,n})^*\left(\mathcal O_{G(p,n)}(1)\right)|_{D^-_j}-(r-1)D^{-j}_{-1}-\cdots-\cdots-D^{-j}_{-(r-1)}\,,\\
\check B^{-j}_r&=D^{-j}_{-r}=(R_{s,p,n})^*\left(\mathcal O_{G(p,n)}(1)\right)|_{D^-_j}-rD^{-j}_{-1}-\cdots-\cdots-2D^{-j}_{-(r-1)}\,.\\
\end{split}
\end{equation} 
Then (\ref{sum1}) and (\ref{sum2}) follow directly.

When $1\leq j\leq r-2$, subtracting (\ref{prdm+1}) from (\ref{prdm}), we derive that
\begin{equation}\label{rdjjv}
D^{-j}_{-1}+D^{-j}_{-2}+\cdots+D^{-j}_{-(j-1)}+D^{-j}_{-j}+D^{-j}_{-(j+1)}=\check B^{-j}_j-\check B^{-j}_{j+1}\,.\\
\end{equation}
Restricting (\ref{tk2}) to $D^-_j$ and plugging in (\ref{rdjjv}), we have that
\begin{equation}
\begin{split}
 -K_{\mathcal T_{s,p,n}}|_{D^{- }_j}=&2\sum_{m=1}^{r-1}\check B^{-j}_m +\sum_{i=1}^{j+1}D^{-j}_{-i}+\sum_{i=j+2}^rD^{-j}_{-i}+\sum_{i=r+2-j}^rD^{-j}_{+i}\,\\
=&2\sum_{m=1}^{r-1}\check B^{-j}_m +\check B^{-j}_j-\check B^{-j}_{j+1}+\sum_{i=j+2}^{r}D^{-j}_{-i}+\sum_{i=r+2-j}^rD^{-j}_{+i}\,.\\
\end{split}
\end{equation}
When $j=r-1$ or $r$, we have that
\begin{equation}
D^{-j}_{-1}+D^{-j}_{-2}+\cdots+D^{-j}_{-(j-1)}+D^{-j}_{-j}+D^{-j}_{-(j+1)}=\check B^{-j}_{r-1};\\
\end{equation}
hence,
\begin{equation}
\begin{split}
 -K_{\mathcal T_{s,p,n}}|_{D^{- }_j}=&2\sum_{m=1}^{r-1}\check B^{-j}_m +\sum_{i=1}^{j+1}D^{-j}_{-i}+\sum_{i=j+2}^rD^{-j}_{-i}+\sum_{i=r+2-j}^rD^{-j}_{+i}\,\\
=&2\sum_{m=1}^{r-1}\check B^{-j}_m +\check B^{-j}_{r-1}+\sum_{i=j+2}^{r}D^{-j}_{-i}+\sum_{i=r+2-j}^rD^{-j}_{+i}\,.\\
\end{split}
\end{equation}

Next we consider the restriction on $D^{+}_j$. By Lemma \ref{bb} and the fact that $D^-_{1},D^-_2,\cdots,D^-_{r+1-j}$ have empty intersection with $D^+_j$, we can conclude that
\begin{equation}\label{pr+dbst}
\begin{split}
\check B^{+j}_0&=D^{+j}_{+r}=(R_{s,p,n})^*\left(\mathcal O_{G(p,n)}(1)\right)|_{D^+_j}-\sum_{i=1}^{r-1}(r+1-i)\cdot D^{+j}_{+i}\,,\\
\check B^{+j}_1&=(R_{s,p,n})^*\left(\mathcal O_{G(p,n)}(1)\right)|_{D^+_j}-\sum_{i=1}^{r-1}(r-i)\cdot D^{+j}_{+i}\,,\\
&\,\,\,\vdots\\
\end{split}
\end{equation}
\begin{equation}\label{pr+dm-2}
\begin{split}
\check B^{+j}_{r-1-j}&=(R_{s,p,n})^*\left(\mathcal O_{G(p,n)}(1)\right)|_{D^+_j}-\sum_{i=1}^{j+1}(j+2-i)\cdot D^{+j}_{+i}\,,\\
\end{split}
\end{equation}
\begin{equation}\label{pr+dm-1}
\begin{split}
\check B^{+j}_{r-j}&=(R_{s,p,n})^*\left(\mathcal O_{G(p,n)}(1)\right)|_{D^+_j}-\sum_{i=1}^{j}(j+1-i)\cdot D^{+j}_{+i}\,,\\
\end{split}
\end{equation}
\begin{equation}\label{pr+dm}
\begin{split}
\check B^{+j}_{r+1-j}&=(R_{s,p,n})^*\left(\mathcal O_{G(p,n)}(1)\right)|_{D^+_j}-\sum_{i=1}^{j-1}(j-i)\cdot D^{+j}_{+i}\,,\\
\end{split}
\end{equation}
\begin{equation}\label{pr+dm+1}
\begin{split}
\check B^{+j}_{r+2-j}&=(R_{s,p,n})^*\left(\mathcal O_{G(p,n)}(1)\right)|_{D^+_j}-\sum_{i=1}^{j-2}(j-1-i)\cdot D^{+j}_{+i}-D^{+j}_{-(r+2-j)}\,,\\
\end{split}
\end{equation}
\begin{equation}\label{p116}
\begin{split}
&\,\,\,\vdots\\
\check B^{+j}_{r-1}&=(R_{s,p,n})^*\left(\mathcal O_{G(p,n)}(1)\right)|_{D^+_j}-D^{+j}_{+1}-\sum_{i=r+2-j}^{r-1}(r-i)\cdot D^{+j}_{-i}\,.\\
\check B^{+j}_r=D^{+j}_{-r}&=(R_{s,p,n})^*\left(\mathcal O_{G(p,n)}(1)\right)|_{D^+_j}-\sum_{i=r+2-j}^{r-1}(r+1-i)\cdot D^{+j}_{-i}\,.\\
\end{split}
\end{equation} 

Then (\ref{pp1}) and (\ref{pp2}) follow from (\ref{pr+dm-2}), (\ref{pr+dm-1}), (\ref{pr+dm}),  (\ref{pr+dm+1}), and (\ref{p116}).

When $1\leq j\leq r-2$, subtracting (\ref{pr+dm-2}) from (\ref{pr+dm-1}), we derive that
\begin{equation}\label{pr+djj}
D^{+j}_{+1}+D^{+j}_{+2}+\cdots+D^{+j}_{+(j-1)}+D^{+j}_{+j}+D^{+j}_{+(j+1)}=\check B^{+j}_{r-j}-\check B^{+j}_{r-1-j}\,.\\
\end{equation}
Restricting (\ref{tk2}) to $D^+_j$ and plugging in (\ref{pr+djj}), we have that
\begin{equation}
\begin{split}
 -K_{\mathcal T_{s,p,n}}|_{D^{+}_j}=&2\sum_{m=1}^{r-1}\check B^{+j}_m +\sum_{i=1}^{j+1}D^{+j}_{+i}+
\sum_{i=j+2}^rD^{+j}_{+i}+\sum_{i=r+2-j}^rD^{+j}_{-i}\,,\\
=&2\sum_{m=1}^{r-1}\check B^{+j}_m +\check B^{+j}_{r-j}-\check B^{+j}_{r-1-j}+\sum_{i=j+2}^{p}D^{+j}_{+i}+\sum_{i=r+2-j}^rD^{+j}_{-i}\,.\\
\end{split}
\end{equation}
When $j=r-1$ or $r$, 
\begin{equation}
D^{+j}_{+1}+D^{+j}_{+2}+\cdots+D^{+j}_{+(j-1)}+D^{+j}_{+j}+D^{+j}_{+(j+1)}=\check B^{+j}_{1};\\
\end{equation}
hence,
\begin{equation}
\begin{split}
 -K_{\mathcal T_{s,p,n}}|_{D^{+ }_j}=&2\sum_{m=1}^{r-1}\check B^{+j}_m +\sum_{i=1}^{j+1}D^{+j}_{+i}+\sum_{i=j+2}^rD^{+j}_{+i}+\sum_{i=r+2-j}^rD^{+j}_{-i}\,\\
=&2\sum_{m=1}^{r-1}\check B^{+j}_m +\check B^{+j}_{1}+\sum_{i=j+2}^{r}D^{+j}_{+i}+\sum_{i=r+2-j}^rD^{+j}_{-i}\,.\\
\end{split}
\end{equation}

We complete the proof of Lemma \ref{rwk3}.\,\,\,\,$\endpf$

\end{document}